\documentclass[11pt,a4paper,USenglish]{article}

\usepackage[utf8]{inputenc}
\usepackage[T1]{fontenc}
\usepackage{lmodern}
\usepackage{microtype}

\usepackage[margin=2.5cm]{geometry}

\usepackage{amsmath}
\usepackage{amssymb}
\usepackage{amsfonts}
\usepackage{amsthm}
\usepackage{mathtools}

\usepackage{graphicx}
\usepackage{float}
\usepackage{booktabs}
\usepackage{caption}
\usepackage{caption}
\usepackage{multirow}
\usepackage{array}

\usepackage[numbers,square,sort&compress]{natbib}

\usepackage[colorlinks=true,
linkcolor=blue,
citecolor=blue,
urlcolor=blue]{hyperref}

\theoremstyle{plain}
\newtheorem{theorem}{Theorem}

\theoremstyle{definition}

\newtheorem{assumption}{Assumption}

\theoremstyle{remark}
\newtheorem{remark}{Remark}

\newcommand{\fntext}[2]{%
	\begingroup
	\renewcommand{\thefootnote}{#1}%
	\footnotetext[0]{#2}%
	\endgroup
}

\title{Efficient iterative methods for the nonlinear compressible Reynolds equation with Navier slip boundary condition}

\author{%
	Salma Raiouadi$^{1,\ast}$ \and
	Mohamed El Alaoui Talibi$^{1,\dagger}$ \and
	Mouna Hilal$^{2,\ddagger}$
}

\date{}

\begin{document}
	\maketitle
	
	\fntext{1}{Department of Mathematics, Laboratory of Mathematics, Modeling and Automatic Systems, Faculty of Sciences Semlalia, Cadi Ayyad University, Marrakech, BP 2390, Morocco.}
	\fntext{2}{Department of Mathematics, Laboratory of Modeling and Combinatorics, Faculty of Applied Sciences, Cadi Ayyad University, Safi, BP 4162, Morocco.}
	\fntext{$\ast$}{E-mail: s.raiouadi.ced@uca.ac.ma}
	\fntext{$\dagger$}{E-mail: elalaoui@uca.ac.ma}
	\fntext{$\ddagger$}{E-mail: mouna.hilal@uca.ac.ma}

\begin{abstract}In this study, we present a comparison of iterative methods for the compressible Reynolds equation with Navier slip boundary conditions, based on a fully compressible cavitation model for lubricated bearings. We establish the existence and uniqueness of the solution to the two-dimensional problem. We introduce the Kirchhoff transformation to reformulate the nonlinear density--pressure relationship into a quasi-linear problem. We develop an efficient numerical approach to solve the two-dimensional nonlinear problem. The methods examined are the Gauss--Seidel method, the Newton--Raphson method and parallelized nonlinear Jacobi method, applied to one- and two-dimensional slider bearing configurations. The overall solution time is shown to depend strongly on the choice of method.

\end{abstract}

\section{Introduction} 

Hydrodynamic lubrication plays a fundamental role in reducing friction and wear in numerous engineering systems such as journal bearings, thrust bearings, seals, and micro-electromechanical devices (MEMS). Under the thin-film approximation, the pressure distribution within the lubricant film is governed by the Reynolds equation, originally derived by \citep{Reynolds1886}, which remains the cornerstone of modern lubrication theory.

The classical Reynolds equation assumes that the lubricant satisfies the no-slip boundary condition at solid surfaces \citep{Bayada1986,Chupin2012}. However, experimental and theoretical studies have shown that partial slip may occur at the fluid--solid interface, particularly on hydrophobic surfaces \citep{Pit1999}, partially wetted lubricated contacts \citep{Spikes2003}, superhydrophobic surfaces \citep{Rothstein2010}, and textured or microstructured surfaces \citep{Lee2014}. Such slip effects can significantly reduce friction and, when properly engineered, improve the load-carrying capacity of lubricated contacts. Among the various slip models, the Navier slip condition, originally introduced by Navier (1823) , is one of the most commonly used boundary conditions for modelling partial slip in lubrication problems. It assumes that the tangential slip velocity is proportional to the wall shear stress through the slip length, which characterizes the fluid--solid interaction \citep{Barrat1999,Vinogradova1999}. Another commonly used slip model is the Tresca slip condition (1864), which assumes that slip occurs only when the tangential shear stress reaches a prescribed critical threshold \citep{Choo2007, Spikes20033, Spikes200333}. In lubrication, heterogeneous slip/no-slip patterns have been widely investigated and shown to improve the hydrodynamic performance of bearings by increasing the load-carrying capacity while reducing friction \citep{Salant2004,Fortier2005,BayadaMeurisse2009}.

Another essential phenomenon in hydrodynamic lubrication is cavitation. In some regions of the lubricating film, the pressure predicted by the classical Reynolds equation may become negative, which is physically unrealistic. To overcome this limitation, several cavitation models have been proposed, such as the Half--Sommerfeld model and the Swift--Stieber model. These models are simple to implement but do not conserve mass and may produce unrealistic pressure distributions. So, Jakobsson, Floberg, and Olsson (JFO) (1965) introduced a mass-conserving cavitation model based on a saturation variable. Among its numerical implementations, the Elrod--Adams algorithm \citep{Elrod1981} is the most widely used. More recently, Bayada and Chupin \citep{Bayada2013} reformulated the JFO model as a fully compressible Reynolds equation, in which the density $\rho$ is related to the pressure through a barotropic law. This formulation provides a smooth transition between vapor, mixture, and liquid phases, recovers the JFO/Elrod--Adams model as a limiting case, and avoids explicit tracking of the cavitation boundary by using a single governing equation over the entire domain. Efficient numerical algorithms for this model were subsequently developed in \citep{Bayada2014,Bayada2018}.

The numerical solution of the compressible Reynolds equation remains challenging because of its strong nonlinearity. The pressure-dependent density $\rho(p)$ enters both the nonlinear diffusion coefficient in the elliptic (diffusive) term on the left-hand side and the convective source term on the right-hand side, resulting in a double nonlinear coupling that significantly complicates the numerical treatment of the equation. Classical fixed-point algorithms generally require a large number of iterations, whereas Newton-type methods involve the repeated assembly and factorization of highly nonlinear Jacobian matrices, making them computationally expensive, especially for fine meshes and two-dimensional problems \citep{Wang2011}. These difficulties become even more pronounced when Navier slip boundary condition are incorporated into compressible cavitation models.

In the present work, we consider the fully compressible cavitation model proposed by G. Bayada and L. Chupin \citep{Bayada2013} combined with the Navier slip boundary condition, a coupling that, to the best of our knowledge, has not previously been investigated. Although G. Bayada \citep{Bayada2018} developed a remarkably fast algorithm for the one-dimensional compressible Reynolds equation with Tresca slip model, no corresponding numerical framework currently exists for the two-dimensional case. The double nonlinearity resulting from the compressible density law and the slip boundary condition constitutes the principal numerical challenge addressed in this work. The objective is to develop an efficient numerical approach capable of solving the two-dimensional problem with a significantly reduced computational cost.

The remainder of this paper is organized as follows. Section~2 introduces the two-dimensional compressible Reynolds equation with Navier slip boundary condition together with the fully compressible cavitation model, and establishes the existence and uniqueness of the solution. Section~3 is devoted to the numerical resolution of the two-dimensional problem. A Kirchhoff transformation is first applied to the pressure-dependent density $\rho(p)$ in order to reduce the double nonlinearity of the compressible Reynolds equation. The transformed problem is then discretized by a finite element method with a dual-mesh approximation of the nonlinear density term. The proposed semi-implicit Jacobi scheme is presented, together with the Gauss--Seidel and Newton--Raphson methods used for comparison. Section~4 reports the numerical results. One-dimensional benchmark problems are first considered, including a comparison with the reference algorithms of \citep{Bayada2014,Bayada2018}. Two-dimensional benchmarks on twin-parabolic and textured surfaces are then used to evaluate the accuracy, convergence behavior, and computational efficiency of the methods. Finally, Section~5 summarizes the main conclusions.

\section{Mathematical model and analysis}
\subsection{Modified Reynolds equation}

\begin{figure}[!ht]
	\centering
	\includegraphics[width=0.48\textwidth]{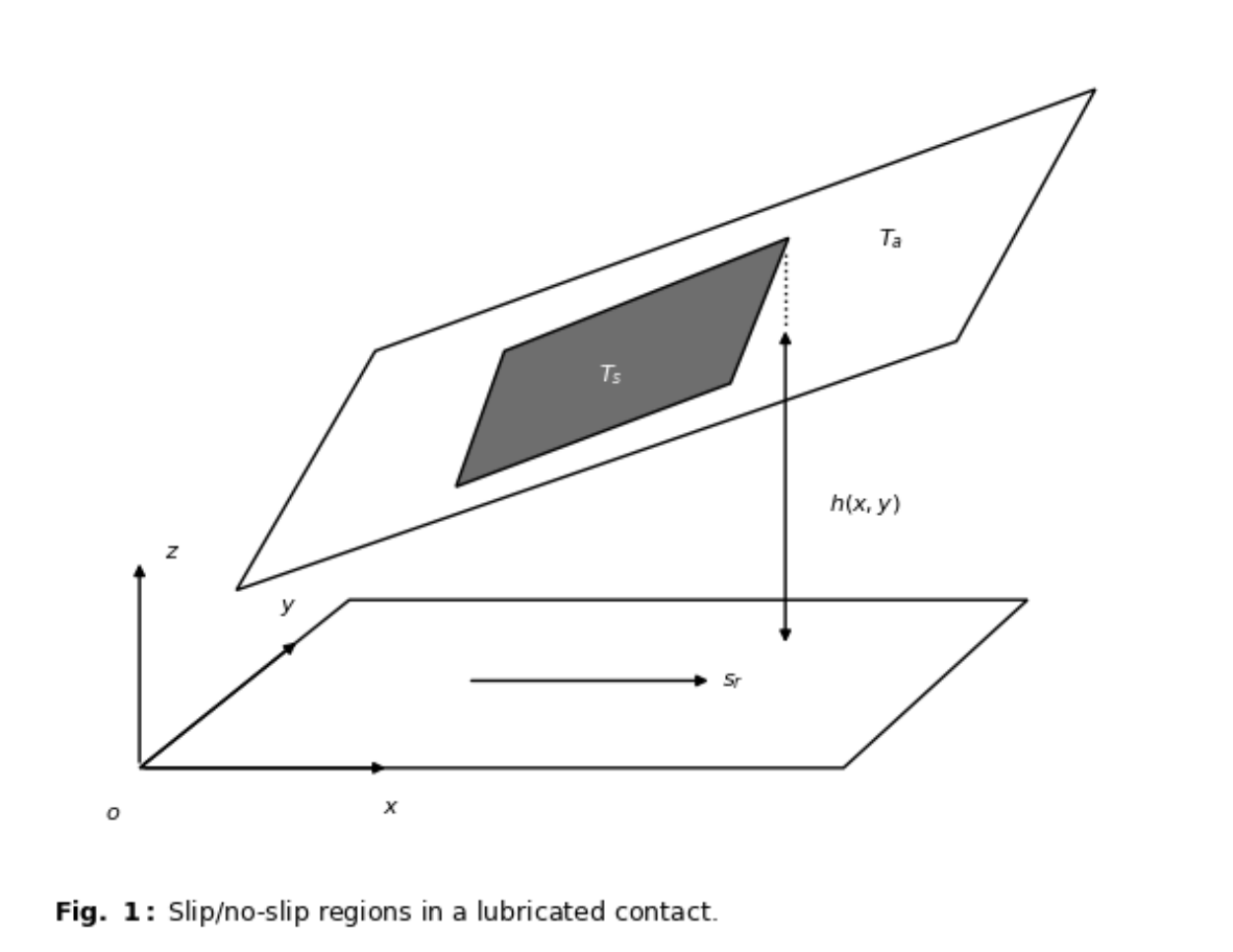}
	\caption{Slip/no-slip regions in a lubricated contact.}
	\label{P1} 
\end{figure}

We consider a lubricated contact shown in Fig.~\ref{P1}, consisting of a stationary top surface and a bottom plane moving in the $x$-direction with velocity $s_r$. The gap between the two surfaces is denoted by $z = h(x,y)$.

The fluid velocity satisfies the following boundary conditions \citep{BayadaMeurisse2009}:
\begin{equation}    
	\begin{cases}
		u(x,y,h(x,y)) = v(x,y,h(x,y)) = w(x,y,h(x,y)) = 0 & \text{on } T_a,\\[0.5em]
		\begin{cases}
			u(x,y,h(x,y)) = -\alpha \mu \dfrac{\partial u}{\partial z}\bigg|_{z=h(x,y)}, \\[0.3em]
			v(x,y,h(x,y)) = -\alpha \mu \dfrac{\partial v}{\partial z}\bigg|_{z=h(x,y)}, \\[0.3em]
			w(x,y,h(x,y)) = 0 
		\end{cases} & \text{on } T_s,\\[0.5em]
		u(x,y,0) = s_r, \quad v(x,y,0) = 0, \quad w(x,y,0) = 0 & \text{on lower surface},
	\end{cases}
	\label{e1}
\end{equation}
where $u$, $v$, and $w$ are the velocity components, $\alpha$ is the slip coefficient, $\mu$ is the lubricant viscosity, and $\alpha\mu$ represents the slip length. The regions $T_a$ and $T_s$ denote the adherence and slip zones on the upper surface, respectively.

Under the thin film approximation, assuming Newtonian fluid, laminar flow, and neglecting external forces, integration of the Stokes equations through the gap yields the modified Reynolds equation \citep{Chupin2012} \citep{BayadaMeurisse2009} : 

\begin{equation}
	\nabla \cdot \left(\left(1+3K(x,y)\right)\rho(p) h^3(x,y)\nabla p \right)=\frac{\partial}{\partial x}\left(6\mu s_r\left(1+K(x,y)\right)\rho(p) h(x,y) \right),
	\label{b8} 
\end{equation}
where the slip coefficient $K(x,y)$ is defined by:
\begin{equation}
	K(x,y)=\begin{cases}
		0 & \text{in no-slip region},\\[0.3em]
		\dfrac{\alpha\mu}{h(x,y)+\alpha\mu} & \text{in slip region}.        
	\end{cases}         
	\label{b9}
\end{equation}

The pressure at the boundary of the domain is equal to the ambient pressure:
\begin{equation}
	p = p_{\text{ext}} \quad \text{on } \partial\Omega,
\end{equation}
where $p_{\text{ext}}$ is a constant.
\subsection{A fully compressible cavitation model}
\label{subsec:cavitation_model}

The cavitation phenomenon is implicitly incorporated into the
density--pressure relation within the framework of vaporous cavitation.
Three distinct regimes are considered: a pure vapor regime, a
liquid--vapor mixture regime, and a pure liquid regime. The input data
are the sound velocities $c_v$ and $c_l$, the vapor and liquid
densities $\rho_v$ and $\rho_l$, and the corresponding viscosities
$\mu_v$ and $\mu_l$ in the two pure regimes. From these parameters, the
transition pressures between the three regimes can be determined,
namely $p_v$ for the vapor--mixture transition and $p_c$ for the
mixture--liquid transition \citep{Bayada2013} \citep{Bayada2014}.

The vapor--mixture transition pressure is given by
\begin{equation}
	p_v = \rho_v c_v^2,
	\label{eq:pv}
\end{equation}
while the mixture--liquid transition pressure is
\begin{equation}
	p_c
	=
	p_v
	+
	N\log\left(
	\frac{\rho_v^2 c_v^2}
	{\rho_l^2 c_l^2}
	\right),
	\label{eq:pc}
\end{equation}
where
\begin{equation}
	N
	=
	\frac{
		\rho_v c_v^2\,\rho_l c_l^2
		(\rho_v-\rho_l)
	}{
		\rho_v^2c_v^2-\rho_l^2c_l^2
	}.
	\label{eq:N}
\end{equation}

It is convenient to introduce the void fraction $\alpha$, defined by
\begin{equation}
	\alpha
	=
	\frac{\rho_l-\rho}
	{\rho_l-\rho_v}.
	\label{eq:alpha}
\end{equation}

The pressure--density relation is then defined by three branches:
\begin{equation}
	p(\rho)
	=
	\begin{cases}
		\rho c_v^2,
		& \rho < \rho_v,\\[0.5em]
		p_c
		+
		N\log\left(
		\frac{
			\rho_v c_v^2\,\rho
		}{
			\rho_l
			\left[
			\rho_v c_v^2(1-\alpha)
			+
			\rho_l c_l^2\alpha
			\right]
		}
		\right),
		& \rho_v \leq \rho \leq \rho_l,\\[0.8em]
		p_c
		+
		c_l^2(\rho-\rho_l),
		& \rho > \rho_l.
	\end{cases}
	\label{eq:fc_pressure_alpha}
\end{equation}

The first and third regimes correspond to the pure vapor and pure liquid regimes, respectively, while the middle regime describes the liquid--vapor mixture. This pressure--density relation provides a continuous description of the transition between the three regimes and allows the cavitation phenomenon to be incorporated into the Reynolds equation without explicitly tracking the cavitation boundary.

\begin{figure}[!ht] \centering \includegraphics[width=0.48\textwidth]{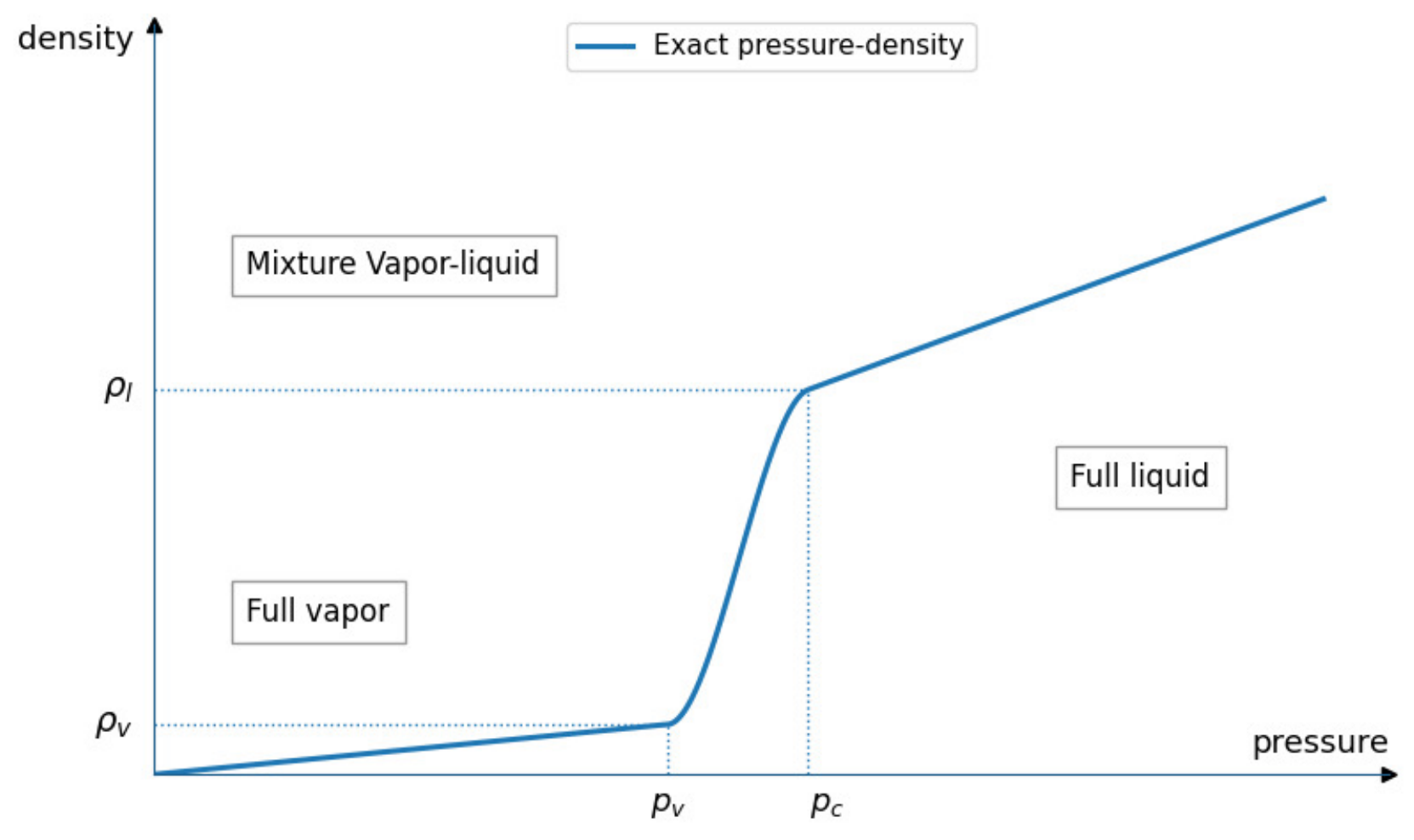} \caption{Density--pressure relation $\rho(p)$ of the fully compressible cavitation model \citep{Bayada2013}. The dashed lines indicate the vapor--mixture ($p_v$) and mixture--liquid ($p_c$) transition pressures.} \label{fig:rho_p} \end{figure}

\subsection{Existence and uniquness results}
\label{subsec:existence}
We consider the compressible Reynolds problem
\begin{equation}
	\begin{cases}
		\text{Find } p \in H^1(\Omega) \text{ such that:}\\[1mm]
		\nabla\cdot\left[
		\left(1+3K(X)\right)h^3(X)\rho(p)\nabla p
		\right]
		=
		\dfrac{\partial}{\partial x}
		\left[
		6\mu s_r\left(1+K(X)\right)h(X)\rho(p)
		\right]
		\quad \text{in }\Omega,\\[1mm]
		p=p_{\mathrm{ext}}
		\quad\text{on }\partial\Omega,
	\end{cases}
	\label{b0}
\end{equation}
where $X=(x,y)$ and $\Omega\subset\mathbb{R}^2$ is a bounded regular domain, under the following assumptions:
\begin{assumption}
	\label{h1}
	The density function $\rho:[0,+\infty[\to[0,+\infty[$ satisfies:
	\begin{itemize}
		\item $\rho$ is continuous and increasing;
		\item $\rho(p)>0$ for all $p>0$;
		\item $\rho(0)=0$.
	\end{itemize}
\end{assumption}
\begin{assumption}
	\label{h2}
	The sliding velocity satisfies $s_r \geq 0$.
\end{assumption}
\begin{assumption}
	\label{h3}
	The film thickness satisfies $h\in L^\infty(\Omega)$ and there exist constants $h_{\min},h_{\max}>0$ such that:
	\[
	0<h_{\min}\leq h(X)\leq h_{\max}<+\infty
	\quad\text{a.e. in }\Omega.
	\]
\end{assumption}

\begin{assumption}
	\label{h4}
	The density function $\rho:[0,+\infty[\to[0,+\infty[$ is Lipschitz continuous.
\end{assumption}

Since $\rho(0)=0$, the diffusion coefficient
$\left(1+3K(X)\right)h^3(X)\rho(p)$ vanishes at $p=0$.
To handle this behavior, we introduce, as in~\cite[Section~4]{buscaglia2005},
the Kirchhoff transformation:

\begin{equation}
	R(s) = \int_0^s \rho(t)\,dt, \qquad s\geq0.
	\label{eq:R_definition}
\end{equation}
By Assumption~\ref{h1}, $R:[0,+\infty[\to[0,+\infty[$ is a continuous, strictly increasing and bijective function. Setting $u=R(p)$, we have
$\nabla u=\rho(p)\nabla p$, and problem~\eqref{b0} is transformed into
the strictly elliptic problem
\begin{equation}
	\begin{cases}
		\nabla\cdot\left[(1+3K(X))h^3(X)\nabla u\right]
		= \dfrac{\partial}{\partial x}\left[ 6\mu s_r(1+K(X))h(X)\,\rho\big(R^{-1}(u)\big)\right]
		&\text{in }\Omega,\\[1mm]
		u=R(p_{\mathrm{ext}}) &\text{on }\partial\Omega,
	\end{cases}
	\label{eq:transformed_problem}
\end{equation}

\begin{theorem}[Existence]
	\label{thm:existence}
	Let Assumptions~\ref{h1}--\ref{h3} hold,
	and assume in addition that
	\begin{equation}
		\rho_l c_l^2 > p_c,
		\label{condcav}
	\end{equation}

	then problem~\eqref{b0} admits at least one non-negative solution.
\end{theorem}

\begin{proof}
	The proof relies on applying theorem~4.1 of Buscaglia \textit{et al.}~\citep{buscaglia2005} to the transformed problem~\eqref{eq:transformed_problem}, obtained through the Kirchhoff transformation~\eqref{eq:R_definition}.
	
	We verify first that the coefficients of~\eqref{eq:transformed_problem} satisfy hypotheses (H1), (H2), and (3.1) of theorem~4.1 in~\citep{buscaglia2005}:
	\begin{itemize}
		\item Hypothesis (H1): Since $h, K \in L^\infty(\Omega)$, $(1+3K(X))h^3(X)$ and $6\mu s_r(1+K(X))h(X)\,\rho\big(R^{-1}(u)\big)$ are measurable with respect to $X$. The map $u \mapsto \rho(R^{-1}(u))$ is continuous on $]0,+\infty[$ as a composition of continuous functions, ensuring $u \mapsto 6\mu s_r(1+K(X))h(X)\,\rho\big(R^{-1}(u)\big)$ is continuous a.e. in $\Omega$.\\
		
		\item Hypothesis (H2): Since \(K(X)\geq0\) and \(h_{\min}\leq h(X)\leq h_{\max}\), the coefficient \((1+3K)h^3\) is bounded and satisfies $h_{\min}^3 \;\leq\; (1+3K(X))h^3(X) \;\leq\; (1+3\|K\|_{L^\infty(\Omega)})\,h_{\max}^3
		\quad\text{a.e. in }\Omega.
		$
		Moreover, for any compact set $\mathcal{K} \subset ]0,+\infty[$, the continuity of $R^{-1}$ and $\rho$ ensures that $C_{\mathcal{K}} := \sup_{u\in\mathcal{K}} \rho(R^{-1}(u)) < +\infty$. Hence, $\sup_{u\in\mathcal{K}} \sup_{X\in\Omega} |6\mu s_r(1+K(X))h(X)\,\rho\big(R^{-1}(u)\big)| \leq 6\mu s_r \|(1+K)h\|_{L^\infty(\Omega)} C_{\mathcal{K}} < +\infty$.
		
		\item Since $\rho(0)=0$ (Assumption~\ref{h1}), hypothesis (3.1) is immediately satisfied.
	\end{itemize}
	
	Theorem~4.1 of~\citep{buscaglia2005} further requires the existence of
	$q>p_{\mathrm{ext}}$ and $\eta>1$ such that
	\begin{equation}
		K_\Omega\,
		\frac{6\mu s_r \|(1+K)h\|_{L^\infty(\Omega)}\,\rho(\eta q)}
		{h_{\min}^3\,\big(R(\eta q)-R(q)\big)}
		\;\leq\;1,
		\label{eq:final_growth}
	\end{equation}
	where $K_\Omega>0$ is a constant depending only on $\Omega$, arising
	from the Sobolev-type estimate in~\citep[Lemma~2.1]{buscaglia2005}.
	
	Since $\rho$ is increasing (Assumption~\ref{h1}), we have, for
	$t\in[q,\eta q]$, $\rho(t)\ge\rho(q)$, so that
	\begin{equation}
		R(\eta q)-R(q)=\int_{q}^{\eta q}\rho(t)\,dt
		\;\geq\;(\eta-1)\,q\,\rho(q).
		\label{cond}
	\end{equation}
	Consequently, it suffices to show that:
	\begin{equation}
		K_\Omega\,
		\frac{6\mu s_r \|(1+K)h\|_{L^\infty(\Omega)}\,\rho(\eta q)}
		{h_{\min}^3\,(\eta-1)\,q\,\rho(q)}\;\leq\;1
		\label{eq:reduced}
	\end{equation}
	to ensure~\eqref{eq:final_growth}.
	
	We now show that condition~\eqref{eq:reduced} can always be satisfied
	for our model. Recall that, for $p>p_c$,
	$\rho(p)=kp+l$ with $k=1/c_l^2>0$ and $l=\rho_l-p_c/c_l^2$.
	
	Let $q>\max(p_{\mathrm{ext}},p_c)$ be large enough, and set:
	\begin{equation}
		\label{eq:gamma_star}
		\eta := \frac{1+B(q)}{1-A(q)},
	\end{equation}
	where
	\begin{equation}
		A(q)=\frac{C_0\,k}{kq+l},
		\qquad
		B(q)=\frac{C_0\,l}{q(kq+l)},
		\label{eq:AB_qstar}
	\end{equation}
	with
	\begin{equation}
		\label{eq:C0}
		C_0:=\frac{K_\Omega\,6\mu s_r\|(1+K)h\|_{L^\infty(\Omega)}}{h_{\min}^3}.
	\end{equation}
	Since $k,l,C_0,q>0$, we have $A(q)>0$ and $A(q)<1$. Moreover, by condition~\eqref{condcav}, i.e.\ $l=\rho_l-p_c/c_l^2>0$, we obtain $B(q)>0$. So that
	\[
	\eta := \frac{1+B(q)}{1-A(q)} \;>\;1.
	\]
	
	Since $\eta q>q>p_c$, both $q$ and $\eta q$ lie in the pure liquid
	region, so $\rho(q)=kq+l$ and $\rho(\eta q)=k\eta q+l$. Substituting
	these into~\eqref{eq:reduced} gives:
	
	\begin{align*}
		C_0\,\frac{\rho(\eta q)}{(\eta-1)\,q\,\rho(q)}
		&= C_0\,\frac{k\eta q + l}{(\eta-1)\,q\,(kq+l)}\\[2mm]
		&= \frac{A(q)\,\eta + B(q)}{\eta - 1},
	\end{align*}
	
	By definition of $\eta$, we obtain:
	\[
	\frac{A(q)\,\eta + B(q)}{\eta - 1} = 1 \;\leq\; 1,
	\]
	so, \eqref{eq:reduced} is verified, and therefore \eqref{eq:final_growth} holds.
	
	Hence, by Theorem~4.1 in~\citep{buscaglia2005}, problem~\eqref{b0} admits at least one positive solution.
\end{proof}

\begin{remark}
	In most physical applications, the numerical values of the parameters imply that condition~\ref{condcav} is satisfied.
\end{remark}

\begin{theorem}[Uniquness]
	\label{thm:uniqueness}
	Under assumptions~\ref{h1}--\ref{h4}, there exists a unique positive
	solution $p \in H^1(\Omega)$ to problem~\eqref{b0}.
\end{theorem}
\begin{proof}
	Under Assumptions~\ref{h1}--\ref{h4}, the hypotheses of Theorem~2.6 in~\citep{buscaglia2005} are satisfied. Hence, the theorem applies directly.
\end{proof}
\section{Numerical methods for the two-dimensional compressible Reynolds problem}
The strong nonlinearity in problem \eqref{b0} makes classical iterative
approaches based on the pressure variable, such as sequential methods
(e.g., Gauss--Seidel iterations), exhibit slow convergence, while
Newton-type methods require the computation of a highly nonlinear Jacobian
matrix. To overcome these difficulties and simplify the structure of \eqref{b0},
we apply the Kirchhoff transformation defined by \eqref{eq:R_definition}. For the numerical resolution, we use the piecewise-linear approximation
of the density--pressure relation proposed by G. Bayada \cite{Bayada2014},
which preserves the overall shape of the fully compressible law while
admitting an explicit, closed-form expression:
\begin{equation}
	\rho(p) = 
	\begin{cases} 
		\rho_{v} & \text{if } p \le p_{v}, \\[0.5em]    
		a p + b & \text{if } p_{v} < p \le p_{c}, \\[0.5em]
		\rho_l & \text{if } p > p_{c},
	\end{cases}
	\label{eq:rho_approx}
\end{equation}

where $\rho_l$ is the liquid density, $p_c$ is the saturation pressure, and $p_v$ is the vapor pressure below which the density is assumed constant and equal to $\rho_v$. Setting

\begin{equation}
	a := \dfrac{\rho_l - \rho_{v}}{p_{c} - p_{v}}
	\label{a}
\end{equation}

\begin{equation}
	b := \dfrac{\rho_v p_{c} - \rho_{l} p_{v}}{p_{c} - p_{v}}.
	\label{b}
\end{equation}

This approximation replaces the exact logarithmic transition in
Eq.~\eqref{eq:fc_pressure_alpha} by a linear one, allowing the density
to be evaluated explicitly without iterative procedures such as Newton's
method or dichotomy. Consequently, for the piecewise-linear density law
\eqref{eq:rho_approx}, both the Kirchhoff transformation $R$, defined
in Eq.~\eqref{eq:R_definition}, and its inverse $R^{-1}$ admit explicit
analytical expressions.

Since $\rho(p)>0$ for all $p\in ]0,+\infty[$, the map $R$ is strictly increasing 
and therefore invertible. For the piecewise linear density 
law \eqref{eq:rho_approx}, both $R$ and $R^{-1}$ are available 
analytically:
\begin{equation}
	R(p) =
	\begin{cases}
		\rho_v\, p, 
		& p \leq p_v, \\[4pt]
		\rho_v p_v + a(p-p_v) + \dfrac{a}{2}(p^2-p_v^2),
		& p_v < p < p_c, \\[4pt]
		R(p_c) + \rho_l(p-p_c),
		& p \geq p_c,
	\end{cases}
	\label{eq:R_explicit}
\end{equation}
\begin{equation}
	R^{-1}(u) =
	\begin{cases}
		\dfrac{u}{\rho_v}, 
		& u \leq R(p_v), \\[6pt]
		\dfrac{-b + \sqrt{b^2 - 2a\bigl(\rho_v p_v 
				- u - \tfrac{a}{2}p_v^2 - bp_v\bigr)}}{a},
		& R(p_v) < u < R(p_c), \\[6pt]
		p_c + \dfrac{u - R(p_c)}{\rho_l},
		& u \geq R(p_c).
	\end{cases}
	\label{eq:Rinverse_explicit}
\end{equation}

The resulting expression for $\rho\bigl(R^{-1}(u)\bigr)$ is:
\begin{equation}
	\rho\bigl(R^{-1}(u)\bigr) =
	\begin{cases}
		\rho_v, 
		& u \leq R(p_v), \\[6pt]
		\sqrt{\rho_{v}^{2} + 2a\bigl(u - \rho_v p_v \bigr)},
		& R(p_v) < u < R(p_c), \\[6pt]
		\rho_l,
		& u \geq R(p_c).
	\end{cases}
	\label{eq:roRinverse_explicit}
\end{equation}

Setting $u_v := R(p_v)$ and $u_c := R(p_c)$

\begin{remark}
	For $u_v<u<u_c$, the quantity under the square root in
	\eqref{eq:Rinverse_explicit},
	\[
	\Delta := b^2 - 2a\Bigl(\rho_v p_v - u - \tfrac{a}{2}p_v^2 - b\,p_v\Bigr),
	\]
	can be rewritten, using $u_v=\rho_v p_v$, as
	\[
	\Delta = (a p_v + b)^2 + 2a\,(u-u_v) = \rho_v^2 + 2a\,(u-u_v),
	\]
	since $a p_v+b=\rho_v$. As $a>0$ and $u>u_v$, this gives
	$\Delta > \rho_v^2 > 0$, so the square root in
	\eqref{eq:Rinverse_explicit} is well defined and positive on $(u_v,u_c)$.
\end{remark}

For the particular density law introduced in~\ref{eq:rho_approx}, we have:
\begin{theorem}
	\label{thm:th6}
	Let assumptions~\ref{h2} and~\ref{h3} hold. Then problem~\eqref{b0} admits at least one
	non-negative solution $p\in H^1(\Omega)$. Moreover, uniqueness is obtained if
	assumption~\ref{h4} is also satisfied.
\end{theorem}

\begin{proof}
	Let $q > \max(p_{\mathrm{ext}}, p_c)$ be large enough, let $\eta \ge C_0/q + 1$; then all states $s \in [q, \eta q]$ lie in the pure liquid regime, where $\rho(s) = \rho_l$. Thus
	$\sup_{s\in[q,\eta q]} \rho(s) = \rho_l$ and $R(\eta q)-R(q) = \rho_l(\eta-1)q$,
	so that condition~\eqref{eq:reduced} reduces to
	$\dfrac{C_0}{(\eta - 1) q} \le 1$, where $C_0$ is defined in~\eqref{eq:C0}.
	This inequality is always satisfied, independently of
	condition~\eqref{condcav}. This proves existence.
	
	Uniqueness then follows from Theorem~\ref{thm:uniqueness}, since
	Assumptions~\ref{h2}--\ref{h4} hold.
\end{proof}

\subsection{Finite element discretization and dual mesh treatment of nonlinear term}
\label{sectionFEM}
To avoid confusion between the film thickness $h(X)$ and the mesh size $h$, we denote the hysical film thickness $h$ and the mesh size $h$, we denote the physical film thickness by $H(X)$ throughout the rest of the paper. The notation $h$ is reserved for the mesh size.

We consider a rectangular computational domain $\Omega=[0,L_x]\times[0,L_y]$,
where $L_x$ and $L_y$ are the characteristic lengths in the $x$- and
$y$-directions, respectively. Introducing the dimensionless variables
\begin{equation}
	\bar x=\frac{x}{L_x},\quad \bar y=\frac{y}{L_y},\quad \bar H=\frac{H}{h_0},\quad
	\bar p=\frac{h_0^2}{6\mu s_r L_x}p,\quad \bar p_{\mathrm{ext}}=\frac{h_0^2}{6\mu s_r L_x}p_{\mathrm{ext}},
\end{equation}
the domain becomes $\hat\Omega=[0,1]\times[0,1]$, and Eq.~\eqref{b0}
takes the dimensionless form

\begin{equation}
	\begin{cases}
		\dfrac{\partial}{\partial \bar{x}}\left(
		(1+3\bar{K})\bar{H}^{3}\bar{\rho}(\bar{p})
		\dfrac{\partial \bar{p}}{\partial \bar{x}}
		\right) + \left(\dfrac{L_x}{L_y}\right)^2 \dfrac{\partial}{\partial \bar{y}}\left(
		(1+3\bar{K})\bar{H}^{3}\bar{\rho}(\bar{p})
		\dfrac{\partial \bar{p}}{\partial \bar{y}}
		\right) \\[0.3cm]
		\qquad = \dfrac{\partial}{\partial \bar{x}}\left(
		(1+\bar{K})\bar{H}\,\bar{\rho}(\bar{p})
		\right) \quad \text{in } \hat{\Omega}, \\[0.3cm]
		\bar{p} = \bar{p}_{\mathrm{ext}} \quad \text{on } \partial\hat{\Omega},
	\end{cases}
	\label{eq:ReynoldsND_full}
\end{equation}
where $\bar{\rho}(\bar{p}):= \rho\left(\dfrac{h_0^2}{6\mu s_r L_x}\bar{p}\right)$ and $\bar K(\bar x,\bar y)=K(L_x\bar x,\;L_y\bar y)$. \\

Applying the Kirchhoff transformation $\tilde{u} = R(\bar{p})$ and setting $\bar{u} = \tilde{u} - R(\bar{p}_{\mathrm{ext}})$, we obtain:
\begin{equation}
	\begin{cases}
		\dfrac{\partial}{\partial\bar{x}}\left[
		(1+3\bar{K})\bar{H}^3
		\dfrac{\partial \bar{u}}{\partial\bar{x}}
		\right]
		+ \left(\dfrac{L_x}{L_y}\right)^{2}
		\dfrac{\partial}{\partial\bar{y}}\left[
		(1+3\bar{K})\bar{H}^3
		\dfrac{\partial \bar{u}}{\partial\bar{y}}
		\right] \\[0.4cm]
		\qquad = \dfrac{\partial}{\partial\bar{x}}\left[
		(1+\bar{K})\bar{H}\,
		\bar{\rho}\left(R^{-1}(
		\bar{u}+R(\bar{p}_{\mathrm{ext}}))\right)
		\right], \quad \text{in } \hat{\Omega}, \\[0.4cm]
		\bar{u} = 0 \quad \text{on } \partial\hat{\Omega}.
	\end{cases}
	\label{eq:shifted_u_2D}
\end{equation}

We define
\begin{align}
	\mathcal{D}(\bar{X}) &= (1+3\bar{K}(\bar{X}))\bar{H}^3(\bar{X}), \\
	F(\bar{X},\bar{u}) &= (1+\bar{K}(\bar{X}))\bar{H}(\bar{X})\,
	\bar{\rho}\left(R^{-1}(
	\bar{u}+R(\bar{p}_{\mathrm{ext}}))\right).
\end{align}

\begin{remark}
	Under the Kirchhoff transformation, the principal part of the operator, corresponding to the diffusion term, becomes linear in $\bar{u}$, with a coefficient $\mathcal{D}(X)$ depending only on the geometry and the slip conditions. The right-hand  side $F(\bar{X},\bar{u})$ remains nonlinear through the density $\bar{\rho}\big(R^{-1}(\bar{u}+R(\bar{p}_{\mathrm{ext}}))\big)$.
\end{remark}

Let $\tau_h = \{ T_k \}_{k=1}^{n_t}$ be a regular triangulation of $\hat{\Omega}$ with mesh size $h > 0$. We define the finite element spaces:
\begin{align}
	H_h &= \{ v_h \in C^0(\hat{\Omega}) : v_h|_{T_k} \in \mathbb{P}^1(T_k), \; \forall T_k \in \tau_h \}, \\
	V_h &= \{ v_h \in H_h : v_h = 0 \text{ on } \partial\hat{\Omega} \}.
\end{align}

Let $\{\phi_i\}_{i=1}^{n_v}$ be the standard nodal basis of $V_h$, where $n_v$ is the number of interior nodes. The discrete problem is to find $\bar{u}_h \in V_h$ such that
\begin{equation}
	\int_{\hat{\Omega}} \mathcal{D} \cdot \dfrac{\partial \bar{u}_h}{\partial \bar{x}} \cdot \dfrac{\partial v_h}{\partial \bar{x}}  \,d\bar{X} + \left(\dfrac{L_x}{L_y}\right)^{2} \int_{\hat{\Omega}} \mathcal{D} \cdot \dfrac{\partial \bar{u}_h}{\partial \bar{y}} \cdot \dfrac{\partial v_h}{\partial \bar{y}}  \,d\bar{X} = \int_{\hat{\Omega}} F(U_h) \cdot \dfrac{\partial v_h}{\partial \bar{x}} \,d\bar{X}, \quad \forall v_h \in V_h,
	\label{G1}
\end{equation}
where $\mathcal{D} = (1+3\bar{K})\bar{H}^3$ and $F(\bar{u}_h) = (1+\bar{K})\bar{H}\,\bar{\rho}(R^{-1}(\bar{u}_h+R(\bar{p}_{\mathrm{ext}})))$.

\begin{figure}[!ht]
	\centering
	\includegraphics[width=0.48\linewidth]{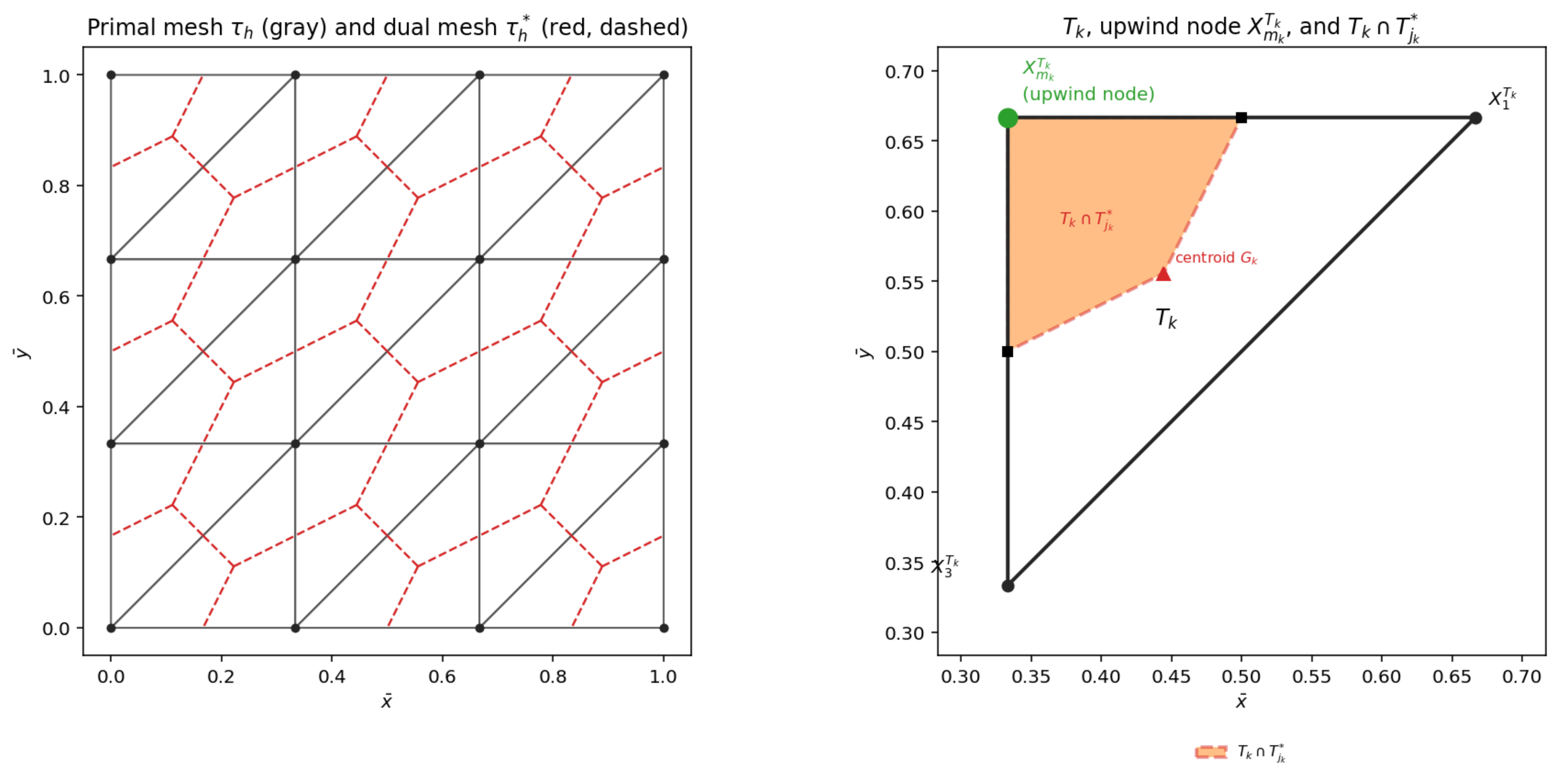}
	\caption{Primal and dual meshes.}
	\label{primaldual}
\end{figure}

The unknown $\bar{u}_h$ is approximated on the primal mesh $\tau_h$ in the usual $\mathbb{P}^1$ sense, while the nonlinear density term is approximated on the dual mesh $\tau_h^* = \{T_j^*\}_{j=1}^{n_v}$, where each cell $T_j^*$ is associated with node $X_j$ (Fig.~\ref{primaldual}):
\begin{equation}
	\bar{u}_h = \sum_{j=1}^{n_v} U_j \phi_j, \qquad \bar{\rho}(R^{-1}(\bar{u}_h+R(\bar{p}_{\mathrm{ext}}))) \approx \sum_{j=1}^{n_v} \rho_j \chi_j,
\end{equation}
where $U_j:=\bar{u}_h(X_j) $ and $\chi_j$ is the characteristic function of the dual cell $T_j^*$ and $\rho_j = \bar{\rho}(R^{-1}(U_j+R(\bar{p}_{\mathrm{ext}})))$. 

Denoting by \(U = \bigl(U_j\bigr)_{1 \leq j \leq n_v} \in \mathbb{R}^n_v\), and substituting into \eqref{G1} and testing with $\phi_i$ yields the following nonlinear discrete system, written in matrix form as
\begin{equation}
	M\,U = E\,\bar{\rho}\big(R^{-1}(U+R(\bar{p}_{\mathrm{ext}}))\big),
	\label{discrete_system}
\end{equation}
where
\[
\bar{\rho}\big(R^{-1}(U+R(\bar{p}_{\mathrm{ext}}))\big)
=
\begin{pmatrix}
	\bar{\rho}\big(R^{-1}(U_1+R(\bar{p}_{\mathrm{ext}}))\big)\\
	\vdots\\
	\bar{\rho}\big(R^{-1}(U_{n_v}+R(\bar{p}_{\mathrm{ext}}))\big)
\end{pmatrix},
\]

\begin{align}
	M_{ij} &= \sum_{T_k \in \tau_h} \int_{T_k} \mathcal{D} \cdot \frac{\partial \phi_i}{\partial \bar{x}} \cdot \frac{\partial \phi_j}{\partial \bar{x}} \, d\bar{X} 
	+ \left(\frac{L_x}{L_y}\right)^{2} \int_{\hat{\Omega}} \mathcal{D} \cdot \frac{\partial \phi_i}{\partial \bar{y}} \cdot \frac{\partial \phi_j}{\partial \bar{y}} \, d\bar{X}, \\
	E_{ij} &= \sum_{T_k \in \tau_h} \int_{T_k \cap T_j^*} (1+\bar{K})\bar{H} \, \frac{\partial \phi_i}{\partial \bar{x}} \, d\bar{X}.
\end{align}

Since $\phi_i$ is $\mathbb{P}^1$, $\nabla\phi_i$ is constant on each triangle $T_k$. Rather than distributing the coefficient $(1+\bar K)\bar H$ among the three dual-cell fractions $T_k\cap T_j^*$, $j\in T_k$, we evaluate it at a single upwind node, defined as the vertex of $T_k$ with the smallest coordinate,

For a triangle $T_k = \operatorname{conv}\{X_1^{T_k}, X_2^{T_k}, X_3^{T_k}\}$, let us write the coordinates of the vertices as $X_m^{T_k} = (\bar{x}_m^{T_k}, \bar{y}_m^{T_k}), \quad m = 1, 2, 3$. 

The upwind node is then defined by $\bar{x}_{m_k}^{T_k} = \min\left\{ \bar{x}_1^{T_k}, \bar{x}_2^{T_k}, \bar{x}_3^{T_k} \right\}$.

Only the dual-cell fraction associated with this upwind node, $T_k \cap T_{j}^*$, contributes, so that for every test function $\phi_i$ of $T_k$
\[
E_{i,j}^{T_k} =
\begin{cases}
	\dfrac{|T_k|}{2} \left(1 + \bar{K}_{m_k}^{T_k}\right) \bar{H}_{m_k}^{T_k} \dfrac{\partial \phi_i}{\partial \bar{x}}, & j = m_k, \\[0.6em]
	0, & j \neq m_k,
\end{cases}
\]
where $|T_k|$ is the area of $T_k$.

\begin{theorem}
	\label{thm:discrete_brouwer}
	Under assumptions~\ref{h1} and~\ref{h3}, the discrete non-linear system~\eqref{discrete_system} admits at least one solution vector $U^\star \in \mathbb{R}^{n_v}$ (corresponding to a discrete solution $\bar{u}_h^\star \in V_h$).
\end{theorem}

\begin{proof}
	We define the bilinear form $a_h(\cdot,\cdot)$ as follows associated with the left-hand side of \eqref{G1}:
	\[
	a_h(\bar{u}_h, v_h) = \int_{\hat{\Omega}} \mathcal{D} \cdot \dfrac{\partial \bar{u}_h}{\partial \bar{x}} \cdot \dfrac{\partial v_h}{\partial \bar{x}}  \,d\bar{X} + \left(\dfrac{L_x}{L_y}\right)^{2} \int_{\hat{\Omega}} \mathcal{D} \cdot \dfrac{\partial \bar{u}_h}{\partial \bar{y}} \cdot \dfrac{\partial v_h}{\partial \bar{y}}  \,d\bar{X}.
	\]
	Since $\mathcal{D}(\bar{X}) = (1+3\bar{K})\bar{H}^3 \ge \bar{H}_{\min}^3 > 0$ a.e.\ in $\hat{\Omega}$, where $\bar{H}_{\min} := \frac{h_{\min}}{h_0}$, Poincaré's inequality guarantees that:
	\begin{equation}
		a_h(v_h, v_h) \ge \bar{H}_{\min}^3 \min\left(1, \left(\frac{L_x}{L_y}\right)^2\right) \|\nabla v_h\|_{L^2(\hat{\Omega})}^2 \ge \dfrac{1}{C_p}\left(\bar{H}_{\min}^3 \min\left(1, \left(\frac{L_x}{L_y}\right)^2\right)\right) \|v_h\|_{L^2(\hat{\Omega})}^2, \quad \forall v_h \in V_h,
		\label{coer}
	\end{equation}
	We set $C_{\mathrm{coer}} = \dfrac{1}{C_p}\left(\bar{H}_{\min}^3 \min\left(1, \left(\frac{L_x}{L_y}\right)^2\right)\right)$, 
	where $C_p>0$ is the Poincaré constant. \\
	Therefore, the bilinear form $a_h(\cdot,\cdot)$ is coercive on $V_h$.
	By coercivity \eqref{coer}, the stiffness matrix $M$ is symmetric positive definite and invertible. System \eqref{discrete_system} reformulates as the fixed-point problem $U = \Phi(U)$, where $\Phi(U) = M^{-1} E \, \bar{\rho}\big(R^{-1}(U + R(\bar{p}_{\mathrm{ext}}))\big)$ is continuous on $\mathbb{R}^{n_v}$. Since $\bar{\rho}(s) \le \rho_l$, we have:
	\[
	\|\Phi(U)\|_2 \le \|M^{-1}\|_2 \, \|E\|_2 \, \rho_l \sqrt{n_v} =: r.
	\]
	Thus, $\Phi$ continuously maps the compact, convex closed ball $\overline{B}(0, r) \subset \mathbb{R}^{n_v}$ into itself. By Brouwer's fixed-point theorem, $\Phi$ admits at least one fixed point $U^\star \in \overline{B}(0, r)$, ensuring the existence of a discrete solution.
\end{proof}

\subsection{Parallelized nonlinear Jacobi method (PNJ)}
We propose a semi-implicit scheme for the nonlinear system
\eqref{discrete_system} that enables a fully vectorized Jacobi
implementation. The diagonal term is evaluated at the current iteration
$m$, while the off-diagonal terms use the values from the previous
iteration $m-1$:
\begin{equation}
	M_{ii}U_i^{(m)} - E_{ii}\,\bar\rho\!\left(R^{-1}\!\left(U_i^{(m)}+R(\bar p_{\mathrm{ext}})\right)\right) = T_{1_i}^{(m-1)},
	\label{eq:jacobi_node}
\end{equation}
where 
\begin{equation}
	T_{1_i}^{(m-1)} = -\sum_{j\ne i} M_{ij}U_j^{(m-1)} + \sum_{j\ne i} E_{ij}\,\bar\rho\!\left(R^{-1}\!\left(U_j^{(m-1)}+R(\bar p_{\mathrm{ext}})\right)\right).
	\label{eq:jacobi_D}
\end{equation}

\begin{theorem}
	\label{thm:existence_uniqueness_jacobi}
	At each node $i \in \{1, \dots, n_v\}$ and for each iteration step $m$, the non-linear equation \eqref{eq:jacobi_node} admits a unique solution $U_i^{(m)} \in \mathbb{R}$.
\end{theorem}
\begin{proof}
	Set $s = U_i^{(m)}$, $T_1 = T_{1_i}^{(m-1)}$, and introduce the function $g : \mathbb{R} \to \mathbb{R}$ defined by:
	\[
	g(s) = M_{ii} \, s - E_{ii} \, \bar{\rho}\left(R^{-1}(s + R_{\mathrm{ext}})\right) - T_1.
	\]
	
	\textbf{Existence:} The function $g$ is continuous on $\mathbb{R}$. Since the density is bounded ($\rho_v \le \bar{\rho} \le \rho_l$) and the diagonal term is strictly positive ($M_{ii} > 0$ by coercivity), it follows that $\lim_{s \to \pm\infty} g(s) = \pm\infty$. The existence of at least one solution $s^* \in \mathbb{R}$ satisfying $g(s^*) = 0$ follows directly from the Intermediate Value Theorem (IVT).
	
	\textbf{Uniqueness:} Since $M_{ii} > 0$, $E_{ii} < 0$ (upwind scheme), and $\bar{\rho}\left(R^{-1}(s + R(\bar{p}_{\mathrm{ext}}))\right)$ is increasing, $g$ is strictly increasing on $\mathbb{R}$, guaranteeing uniqueness.
\end{proof}

We set $u_v := R(p_v)$ and $u_c := R(p_c)$ in dimensionless form: $\bar{u}_v = R(\bar{p}_v)$ and $\bar{u}_c = R(\bar{p}_c)$, where:
\begin{equation}
	\bar{p}_v = \frac{h_0^2}{6\mu s_r L_x} p_v \quad \text{and} \quad \bar{p}_c = \frac{h_0^2}{6\mu s_r L_x} p_c.
\end{equation}
There are three possible cases, depending on the regime considered in Eq.~\eqref{eq:roRinverse_explicit}:
\begin{gather}
	\text{Pure liquid regime } \quad U_{1_i} = \dfrac{T_{1_i}^{m-1}+E_{ii}\,\rho_l}{M_{ii}}, \label{eq:prop1} \\[1ex]
	\text{Pure vapor regime } : \quad U_{2_i} = \dfrac{T_{1_i}^{m-1}+E_{ii}\,\rho_v}{M_{ii}}, \label{eq:prop2} \\[1ex]
	\text{Mixed regime } : \quad U_{3_i} = \dfrac{-\beta_i-\sqrt{\beta_i^2-4\alpha_i\gamma_i}}{2\alpha_i}, \label{eq:prop3}
\end{gather}
where
\begin{align}
	\alpha_i
	&=M_{ii}^{\,2},\label{eq:alpha_i}\\[0.5em]
	\beta_i
	&=-2\Bigl[M_{ii}\,T_{1_i}^{m-1}+a\,E_{ii}^{\,2}\Bigr],\label{eq:beta_i}\\[0.5em]
	\gamma_i
	&=(T_{1_i}^{m-1})^2-E_{ii}^{\,2}\rho_v^2
	-2a\,E_{ii}^{\,2}\Bigl[R(\bar p_{\mathrm{ext}})-u_v\Bigr].\label{eq:gamma_i}
\end{align}

\begin{remark}
	\label{sol+}
	In the mixed regime, Eq.~\eqref{eq:jacobi_node} reads
	\begin{equation}
		M_{ii}U^{(m)}_i-T_{1_i}^{(m-1)}=E_{ii}\sqrt{s(U^{(m)}_i)},\qquad s(U^{(m)}_i)=\rho_v^2+2a\bigl(U^{(m)}_i+R(\bar p_{\mathrm{ext}})-\bar u_v\bigr).
		\label{eq:mixed_sqrt}
	\end{equation}
	Since $E_{ii} < 0$ and $\sqrt{s}\ge 0$, any solution satisfies $U^{(m)}_{i}\le T_{1_i}^{(m-1)}/M_{ii}$.
	Squaring both sides of Eq.~\eqref{eq:mixed_sqrt} yields the quadratic $\alpha_i(U^{(m)}_{i})^{2}+\beta_i U^{(m)}_{i}+\gamma_i=0$, whose roots are
	\[
	U^{(m)}_{i_\pm}=\frac{-\beta_i\pm\sqrt{\beta_i^2-4\alpha_i\gamma_i}}{2\alpha_i}
	=\frac{T_{1_i}^{(m-1)}}{M_{ii}}+\frac{aE_{ii}^{\,2}\pm\tfrac12\sqrt{\beta_i^2-4\alpha_i\gamma_i}}{M_{ii}^{2}},
	\]
	by \eqref{eq:alpha_i} and \eqref{eq:beta_i}. Hence,
	\[
	U^{(m)}_{i_+}-\frac{T_{1_i}^{(m-1)}}{M_{ii}}=\frac{aE_{ii}^{\,2}+\tfrac12\sqrt{\beta_i^2-4\alpha_i\gamma_i}}{M_{ii}^{2}}>0,
	\]
	which is a contradiction. Hence $U^{(m)}_{i_+}$ cannot be the solution in the mixed regime, and the solution in this regime is
	$U_{3_i}=\dfrac{-\beta_i-\sqrt{\beta_i^2-4\alpha_i\gamma_i}}{2\alpha_i}$.
\end{remark}

Once these three possible solutions \ref{eq:prop1}, \ref{eq:prop2}, \ref{eq:prop3} have been computed, it suffices to determine which one is consistent with the corresponding pressure regime:
\begin{equation}
	U_i^m = 
	\begin{cases} 
		U_{1_i}^m, & \text{if } U_{1_i}^m + R(\bar{p}_{\mathrm{ext}}) \le \bar{u}_v, \\[1ex]
		U_{2_i}^m, & \text{if } U_{2_i}^m + R(\bar{p}_{\mathrm{ext}}) \ge \bar{u}_c, \\[1ex]
		U_{3_i}^m, & \bar u_v<U_i^m+R(\bar p_{\mathrm{ext}})<\bar u_c.
	\end{cases} 
	\label{eq:regime_selection}
\end{equation}
Since the matrices $M$ and $E$ are independent of $U$, therefore they are precomputed once. At iteration $m$, the term $T_{1_i}^{(m-1)}$ depends only on the known vector $U^{m-1}$, therefore, it can be computed simultaneously for all nodes:
\begin{equation}
	T_1^{(m-1)} = -M_{\mathrm{off}}\,U^{(m-1)}
	+ E_{\mathrm{off}}\,\bar\rho\!\left(
	R^{-1}\!\left(U^{(m-1)}+R(\bar p_{\mathrm{ext}})\right)
	\right),
	\label{eq:T1_vectorized}
\end{equation}
where $M_{\mathrm{off}}$ and $E_{\mathrm{off}}$ denote the off-diagonal of the matrices $M$ and $E$.

The three candidates $U_k^m$, $k=1,2,3$, are then evaluated for the entire vector at once, with no loop over the nodes: each node is automatically assigned the value $U_{1_i}$, $U_{2_i}$, or $U_{3_i}$ according to its own pressure regime, while the operation itself is performed on the whole vector simultaneously (see algorithm~\ref{alg:jacobi_vectorized} in appendix~\ref{app:algorithms}).
\subsection{Gauss--Seidel scheme for solving the nonlinear system (NGSU)}
\label{sec:GSu}
We solve the nonlinear discrete system \eqref{discrete_system} using a
Gauss--Seidel half-implicit procedure proposed in \citep{Bayada2014},
applied to the Kirchhoff variable $U$. For a forward sweep, i.e.,
$i=1,\ldots,n_v$, the components already updated at iteration $m$ are
those with $j<i$, whereas the components with $j>i$ are taken from the
previous iteration. Thus, the $i$th discrete equation can be written as
\begin{equation}
	M_{ii}U_i^{(m)} - E_{ii}\,\bar\rho\!\left(R^{-1}\!\left(U_i^{(m)}+R(\bar p_{\mathrm{ext}})\right)\right) = T_{2_i}^{(m)},
	\label{eq:GS_u_local}
\end{equation}
where
\begin{align}
	T_{2_i}^{(m)} &= -\sum_{j<i} M_{ij}U_j^{(m)} - \sum_{j>i} M_{ij}U_j^{(m-1)} \notag \\
	&\quad + \sum_{j<i} E_{ij}\,\bar\rho\!\left(R^{-1}\!\left(U_j^{(m)}+R(\bar p_{\mathrm{ext}})\right)\right) \notag \\
	&\quad + \sum_{j>i} E_{ij}\,\bar\rho\!\left(R^{-1}\!\left(U_j^{(m-1)}+R(\bar p_{\mathrm{ext}})\right)\right).
	\label{eq:T2}
\end{align}

\begin{theorem}
	\label{thm:existence_uniqueness_gs}
	At each node $i \in \{1, \dots, n_v\}$ and for each iteration step $m$, the non-linear equation \eqref{eq:GS_u_local} admits a unique solution $U_i^{(m)} \in \mathbb{R}$.
\end{theorem}

\begin{proof}
	Since equation \eqref{eq:GS_u_local} shares the exact same left-hand side operator as \eqref{eq:jacobi_node} with a constant right-hand side $T_{2_i}^{(m)}$, the existence and uniqueness of $U_i^{(m)}$ at each node follow directly from Theorem~\ref{thm:existence_uniqueness_jacobi}.
\end{proof}

Hence, at each node, the nonlinear system
\eqref{eq:GS_u_local} is solved by considering the corresponding
regime of the density--pressure relation.
\begin{gather}
	\text{Pure liquid regime } \quad U_{1_i} = \dfrac{T_{2_i}^{m-1}+E_{ii}\,\rho_l}{M_{ii}}, \label{eq:prop11} \\[1ex]
	\text{Pure vapor regime } : \quad U_{2_i} = \dfrac{T_{2_i}^{m-1}+E_{ii}\,\rho_v}{M_{ii}}, \label{eq:prop22} \\[1ex]
	\text{Mixed regime } : \quad U_{3_i} = \dfrac{-\kappa_i-\sqrt{\kappa_i^2-4\alpha_i\lambda_i}}{2\alpha_i}, \label{eq:prop33}
\end{gather}
where
\begin{align}
	\alpha_i
	&=M_{ii}^{\,2},\\[0.5em]
	\kappa_i
	&=-2\Bigl[M_{ii}\,T_{2_i}^m+a\,E_{ii}^{\,2}\Bigr],\\[0.5em]
	\lambda_i
	&=(T_{2_i}^m)^2-E_{ii}^{\,2}\rho_v^2
	-2a\,E_{ii}^{\,2}\Bigl[R(\bar p_{\mathrm{ext}})-u_v\Bigr].
\end{align}

\begin{remark}
	The argument of Remark~\ref{sol+} also applies to the mixed regime of the Gauss--Seidel scheme, with $T_{1_i}^{(m-1)}$, $\beta_i$, and $\gamma_i$ replaced by $T_{2_i}^{(m)}$, $\kappa_i$, and $\lambda_i$, respectively. Since it relies only on $M_{ii}>0$, $E_{ii} < 0$, and $\sqrt{s}\ge 0$, the root with the $+$ sign is excluded, and the solution is therefore given by Eq.~\eqref{eq:prop33}.
\end{remark}

The three candidates are then tested against their corresponding regime conditions, and the admissible value is selected as the Gauss--Seidel update of $U_i^{(m)}$. To accelerate the convergence of this method, we follow the over-relaxation approach used in~\citep{Bayada2014}. To this end, the solution $U_i^{(m)}$ is not computed directly from $U_{1_i}$ in Eq.~\eqref{eq:prop11}, but is instead given by:

\begin{equation}
	U_i^{(m)} = \omega\,U_{1_i} + (1-\omega)\,U_i^{(m-1)},
\end{equation}
where $\omega \in [1,2[$ is the over-relaxation parameter.
The numerical procedure is summarized in Algorithm~\ref{alg:gs_u} in appendix~\ref{app:algorithms}).

\begin{remark}
	\label{rem:GS_comparison}
	We can also consider the Gauss--Seidel method applied to the dimensionless problem~\eqref{eq:ReynoldsND_full} in the pressure variable $\bar{p}$; we denote this method by NGSP. The corresponding method applied to the Kirchhoff-transformed problem~\eqref{eq:shifted_u_2D} in the variable $\tilde{u}=R(\bar{p})$ (Section~\ref{sec:GSu}) is denoted by NGSU. As for NGSU, the nonlinear diffusion coefficient is handled through a node-by-node nonlinear resolution in $\bar{p}$, following the same principle as in Eqs.~\eqref{eq:prop11}--\eqref{eq:prop33}, but without the benefit of the Kirchhoff transformation. The main computational difference between NGSP and NGSU therefore lies
	in the treatment of the nonlinear diffusion coefficient. In NBSP, the
	matrix associated with the diffusion operator depends on the current
	pressure iterate and must consequently be reassembled at every iteration.
	In NGSU, the Kirchhoff transformation removes this nonlinearity from the
	diffusion operator, so that the matrix $M$ is assembled only once; the
	nonlinear dependence is instead retained solely in the right-hand side,
	through $\bar{\rho}\big(R^{-1}(U+R(\bar{p}_{\mathrm{ext}}))\big)$.
\end{remark}

\subsection{Newton--Raphson method}
\label{subsec:newton}
We also consider the classical Newton--Raphson method applied to the transformed nonlinear system \eqref{discrete_system}.
Let
\[
\mathcal{R}(U)=MU-E\rho(U).
\]
where
$$
M=(M_{ij})_{1\leq i,j\leq n_v},
\qquad
E=(E_{ij})_{1\leq i,j\leq n_v}.
$$
Starting from an initial guess $U^{(0)}$, Newton's method computes
\begin{equation}
	J(U^{(m)})\,\delta U^{(m)}
	=
	-\mathcal{R}(U^{(m)}),
	\qquad
	U^{(m+1)}=U^{(m)}+\delta U^{(m)}.
	\label{eq:newton_step}
\end{equation}
The advantage of the Kirchhoff formulation is that the jacobian matrix
remains particularly simple:
\begin{equation}
	J(U)=M-E\,D_{\rho'}(U),
	\label{eq:newton_jacobian}
\end{equation}
where $D_{\rho'}(U)$ is the diagonal matrix defined by
\[
(D_{\rho'})_{ii}
=
\frac{\rho'(p_i)}{\rho(p_i)}.
\]
To improve the robustness of Newton's method for different slip-zone configurations, we employ an adaptive relaxation approach, as described in Subsection~\ref{subsec:newton_primal_dual}.

The numerical procedure including the adaptive relaxation, is summarized in Algorithm~\ref{alg:newton_adaptive} in appendix~\ref{app:algorithms}).

\section{Numerical resuts}
In this section, we assess the performance of the proposed numerical methods on one- and two-dimensional benchmark problem based on twin-parabolic and textured surface geometries, including textured surfaces. The objectives are to validate the proposed algorithms against reference methods available in the literature and to compare their computational efficiency, convergence behavior, and numerical accuracy.

All numerical experiments were performed on a 64-bit desktop computer equipped with an Intel\textsuperscript{\textregistered} Core\texttrademark~i7-8850H processor running at 2.60~GHz and 16~GB of RAM. The algorithms were implemented in Python using the NumPy and SciPy libraries. 

\subsection{One-dimensional benchmark problem: twin-parabolic slider}
\label{sec:benchmark_1D}
We consider a configuration composed of two successive parabolic sliders, separated by a flat surface. Geometrical and kinematic data are taken from \citep{Bayada2018} and are given in Table~\ref{T2}.

\begin{table}[H]
	\centering
	\caption{Physical twin journal bearing data.}
	\label{T2}
	\renewcommand{\arraystretch}{1.2}
	\begin{tabular}{l l}
		\toprule
		\textbf{Parameter} & \textbf{Value} \\
		\midrule
		Sliding velocity $sr$ (m/s) & 0.3 \\
		Viscosity $\mu$ (Pa s) & 0.04 \\
		Length $L_x$ (m) & 1 \\
		Maximum height $h_{\text{max}}$ (m) & $5 \times 10^{-5}$ \\
		Minimum height $h_{\text{min}_1}$ (m) & $2.5 \times 10^{-5}$ \\
		Minimum height $h_{\text{min}_2}$ (m) & $3 \times 10^{-5}$ \\
		$h_0$ & $h_{\text{max}}$ \\
		Slip coefficient $\alpha$ & 1 \\
		Boundary fluid pressure (Pa) & $10^{5}$ \\
		\bottomrule
	\end{tabular}
\end{table}
\begin{figure}[H]
	\centering
	\includegraphics[width=0.48\textwidth]{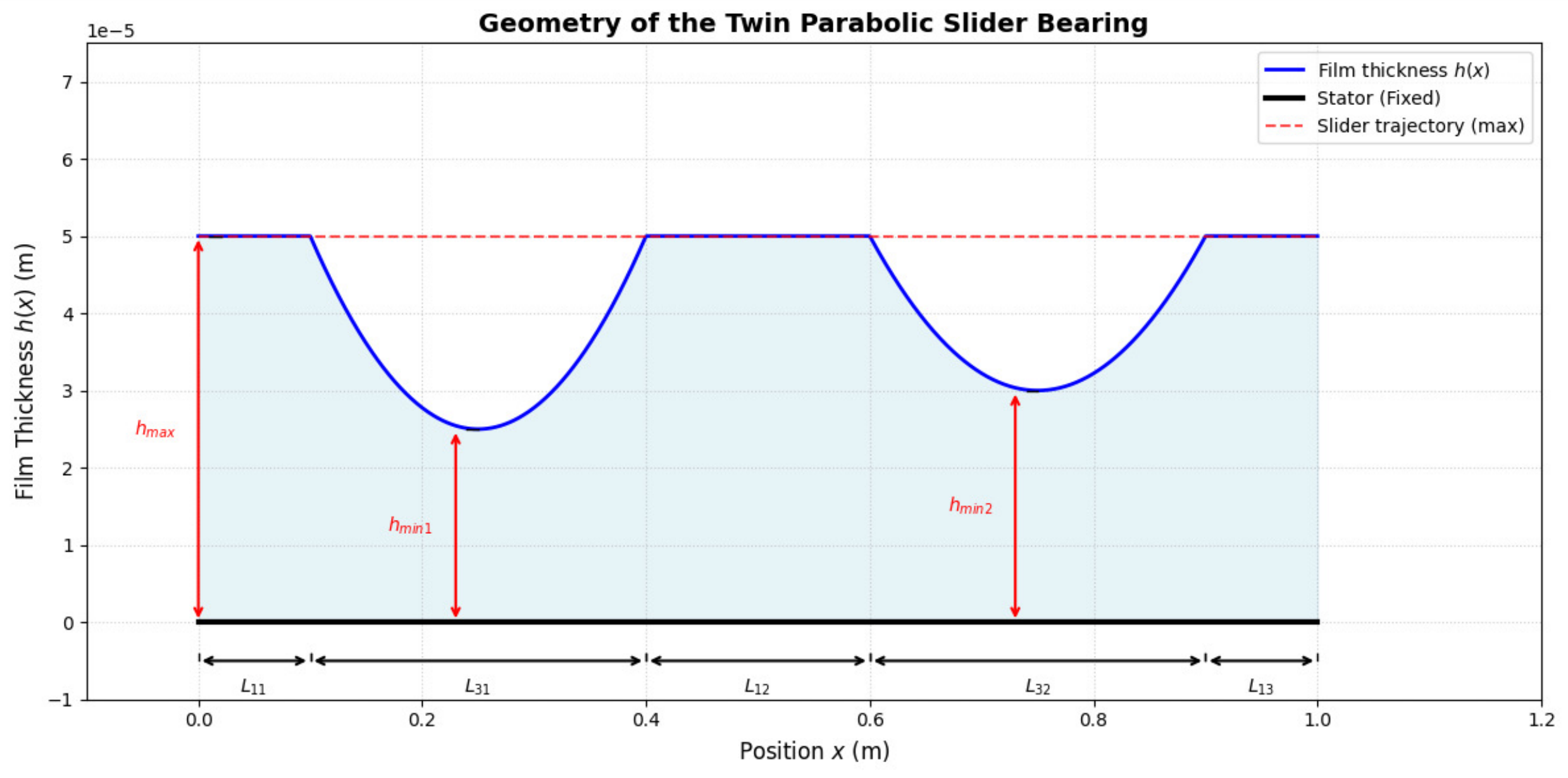}
	\caption{Twin-parabolic slider geometry}
	\label{GTB} 
\end{figure}
The pressure-density relationship is governed by the following parameters: 
\[
\rho_l = 850 \, \text{kg/m}^3, \quad \rho_v = 0.017 \, \text{kg/m}^3, \quad p_c = 82 \times 10^3 \, \text{Pa}, \quad p_v = 47 \times 10^3 \, \text{Pa}
\]

\subsubsection{Reference solvers}

We first validate our implementation of NGSP and NGSU in the one-dimensional
setting, following the semi-implicit Gauss--Seidel algorithm originally proposed by G. Bayada
\citep{Bayada2014} for the compressible Reynolds equation with a classical no-slip
boundary condition. The original finite difference discretization is replaced by a
finite element formulation: the diffusion term is assembled using the trapezoidal rule,
whereas the nonlinear term is approximated by the left rectangle rule instead of the
finite difference approximation used in \citep{Bayada2014}. This approximation preserves
the local dependence of the nonlinear term on the current nodal unknown, and the same
strategy is employed for both NGSP and NGSU.

As a reference benchmark solver, we also consider an algorithm proposed by
\citep{Bayada2018}, hereafter referred to as Fast. The original method
reformulates the one-dimensional Reynolds equation as a sequence of ordinary
differential equations based on the conservation of the mass flow rate, making it
applicable to several slip models, with or without cavitation. In the present work,
this algorithm is adapted to the compressible Reynolds equation with Navier slip
boundary condition through the modified slip coefficients. The unknown flow rate $Q$
is determined from the compatibility condition $F(Q)=0$ associated with the boundary
conditions $p(0)=p(L_x)=p_{\mathrm{ext}}$. Unlike the fixed-point procedure proposed in
Remark~4 of \citep{Bayada2018}, we solve this equation directly by numerical
integration. Finally, the Reynolds ODE is solved backward from $p(L_x)=p_{\mathrm{ext}}$
using an implicit marching procedure. Owing to the piecewise linear density
approximation, the pressure at each node is obtained explicitly in the vapor, mixture,
and liquid regimes, leading to a simple and efficient node-by-node algorithm while
preserving the computational efficiency of the original method.

\subsubsection{Computational performance}
The grid size used throughout this study is $N=500$, with a stopping tolerance of $10^{-8}$. The relaxation parameter $\omega$ is set to $0.1$ for Newton's method and to $1.98$ for the Gauss--Seidel method. The performance of the numerical methods is evaluated for the four slip configurations in terms of the number of iterations, CPU time, and final error, defined as the relative variation $e:=|U^{(m)}-U^{(m-1)}|_2/|U^{(m)}|_2$ and the dimensionless load-carrying capacity $\bar W$. The load-carrying capacity $W$ and its dimensionless form $\bar W$ are defined by:
\begin{equation}
	W=\int_{0}^{L_x}p\,\mathrm{d}x,
	\qquad
	\bar W=\frac{h_0^{2}}{6\mu s_r L_x^{2}}\,W=\int_{0}^{1}\bar p\,\mathrm{d}\bar x,
\end{equation}
The results are reported in Table~\ref{tab:slip_regions1D}, while the corresponding dimensionless pressure profiles are shown in Fig.~\ref{fig:pressure}. In all cases, the final error is of the order of $10^{-8}$, ensuring comparable accuracy among the different methods.

Table~\ref{tab:slip_regions1D} shows that the Kirchhoff transformation provides a clear computational advantage for the Gauss--Seidel method with over-relaxation $ \omega =1.93$. Although NGSP and NGSU require a comparable number of iterations, the CPU time is substantially reduced when the transformed variable is used, with reductions of approximately $52.4\%$, $52.8\%$, $49.8\%$, and $75.0\%$ for the no-slip, full-slip, single-slip, and two-slip configurations, respectively. This indicates that the Kirchhoff transformation mainly reduces the computational cost per iteration rather than the number of iterations.

The parallelized nonlinear Jacobi method requires considerably more iterations than Gauss--Seidel, but its CPU time remains competitive due to the simultaneous nodal updates. For example, in the no-slip case, PNJ requires $40779$ iterations in $5.73$~s, compared with $186.24$~s for NGSP. 

Newton's method provides the lowest CPU time among the considered methods, ranging from approximately $0.13$ to $0.18$~s for the four slip configurations, with a final error of order $10^{-9}$. It is therefore significantly faster than both the Gauss--Seidel and parallelized nonlinear Jacobi method, as well as the Fast algorithm of G. Bayada~\citep{Bayada2018}, whose CPU time ranges from $0.41$ to $4.78$~s. Although the Fast algorithm requires fewer iterations, its computational cost per iteration is higher, particularly for the more complex two-slip configuration.
\begin{table}[H]
	\centering
	\caption{Comparison of the numerical methods for different slip regions (1D case).}
	\label{tab:slip_regions1D}
	\renewcommand{\arraystretch}{1.2}
	\begin{tabular*}{\textwidth}{@{\extracolsep\fill}llcccc}
		\toprule
		Slip region & Method & Iterations & CPU (s) & $\bar W$ & Final error \\
		\midrule
		
		\multirow{5}{*}{No slip}
		& NGSP   & $762$   & $186.2386$ & $0.0139$  & $9.71 \times 10^{-9}$    \\
		& NGSU   & $762$   & $88.6251$  & $0.0139$  & $9.71 \times 10^{-9}$     \\
		& PNJ  & $40779$ & $5.72$     & $0.0139$  & $9.998 \times 10^{-9}$     \\
		& Newton  & $158$   & $0.1349$   & $0.0139$ &  $9.66 \times 10^{-9}$      \\
		& Fast    & $36$    & $1.6994$   &  $0.0140$ & $7.98 \times 10^{-9}$        \\
		\midrule
		
		\multirow{5}{*}{$[0,1]$}
		& NGSP   & $1059$   & $251.9859$  & $0.0083$  & $9.88 \times 10^{-9}$     \\
		& NGSU   & $1042$   & $119.0389$  & $0.0083$  & $9.55 \times 10^{-9}$      \\
		& PNJ  & $41280$  & $5.73$      & $0.0083$  & $9.999 \times 10^{-9}$      \\
		& Newton  & $157$    & $0.1351$    & $0.0083$ & $9.6676 \times 10^{-9}$      \\
		& Fast    & $36$     & $1.5594$    &  $0.0084$ & $7.14 \times 10^{-9}$         \\
		\midrule
		
		\multirow{5}{*}{$[0,0.5]$}
		& NGSP   & $5499$   & $1150.8742$ & $0.0840$  & $1.00 \times 10^{-8}$    \\
		& NGSU   & $5521$   & $577.6943$  & $0.0840$  & $9.99 \times 10^{-9}$     \\
		& PNJ  & $710602$ & $100.09$    & $0.0840$  & $1.00 \times 10^{-8}$      \\
		& Newton  & $155$    & $0.1378$    & $0.0840$ & $9.564 \times 10^{-9}$      \\
		& Fast    & $9$      & $0.4119$    &  $0.0841$ & $5.40 \times 10^{-9}$        \\
		\midrule
		
		\multirow{5}{*}{$[0.18,0.35]\cup[0.68,0.82]$}
		& NGSP   & $1666$   & $765.0220$     & $0.0057$  & $9.99 \times 10^{-9}$    \\
		& NGSU   & $1554$   & $191.4309$     & $0.0057$  & $9.96 \times 10^{-9}$     \\
		& PNJ  & $53033$  & $8.76$         & $0.0057$  & $9.998 \times 10^{-9}$      \\
		& Newton  & $161$    & $0.1805$       & $0.0057$ & $9.36 \times 10^{-9}$      \\
		& Fast    & $81$     & $4.7812$       & $0.0062$ & $8.995 \times 10^{-9}$        \\
		\bottomrule
	\end{tabular*}
\end{table}
\begin{figure}[H]
	\centering
	
	\includegraphics[width=0.48\textwidth]{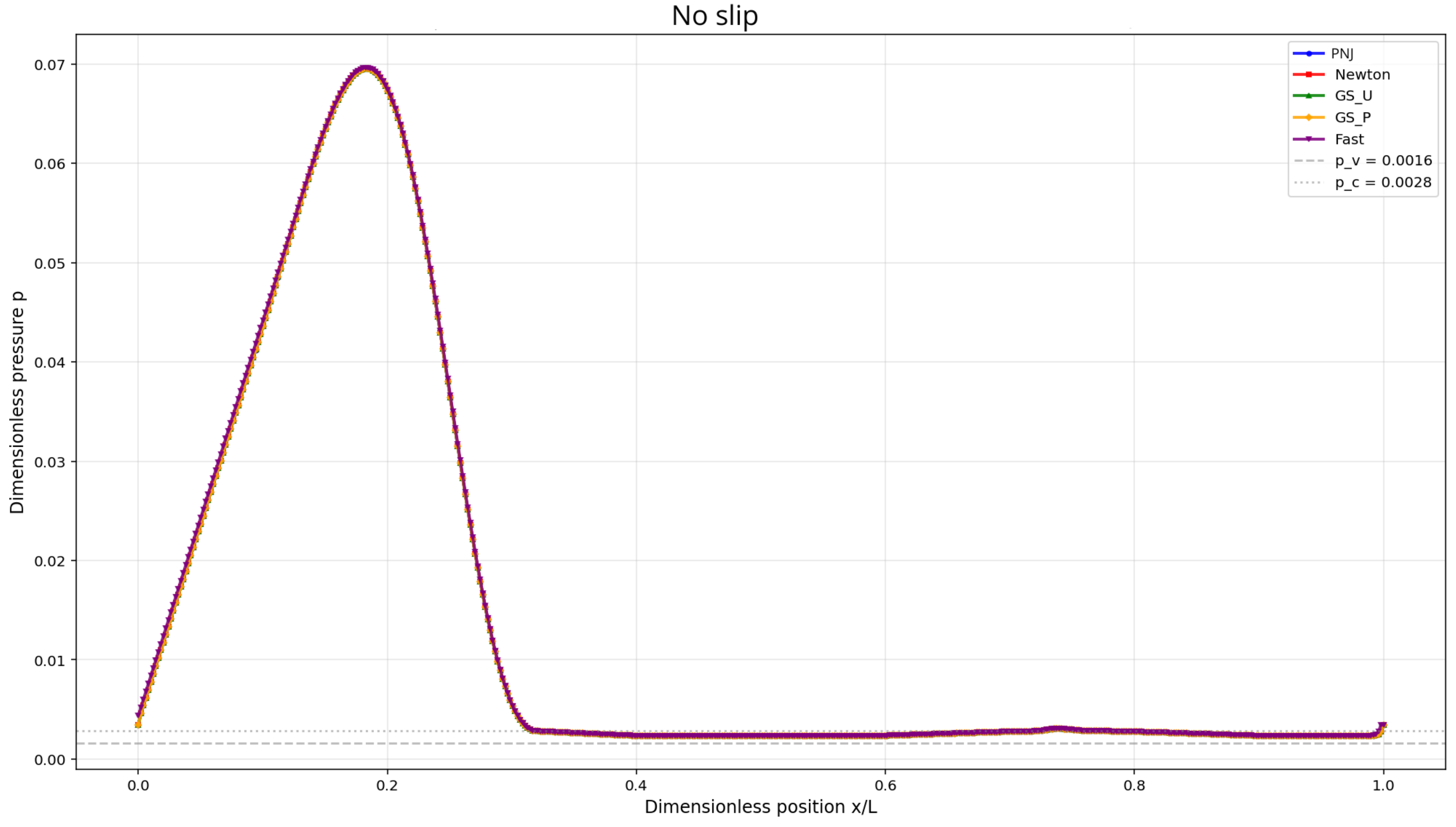}
	\hfill
	\includegraphics[width=0.48\textwidth]{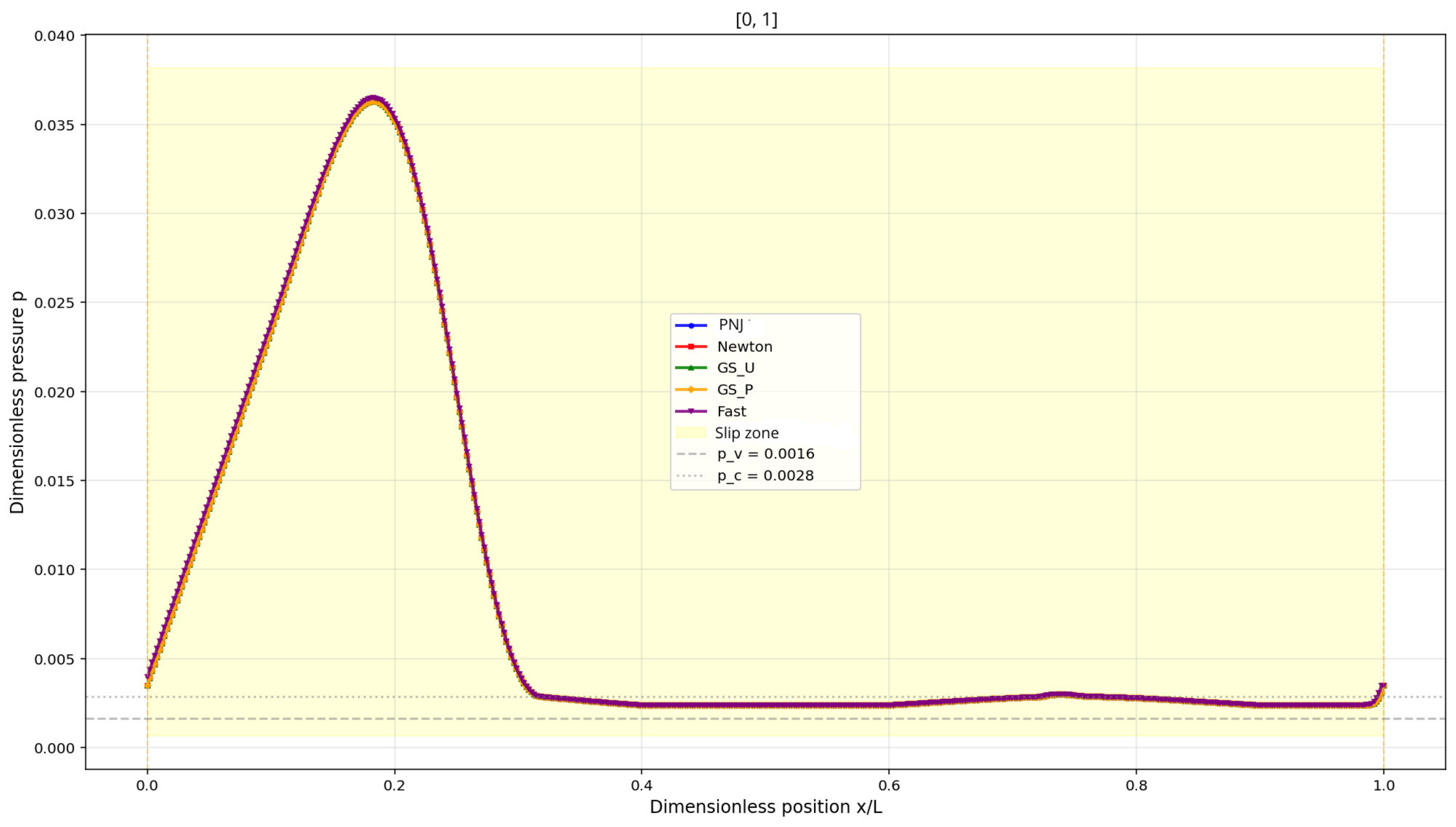}
	
	\vspace{0.4cm}
	
	\includegraphics[width=0.48\textwidth]{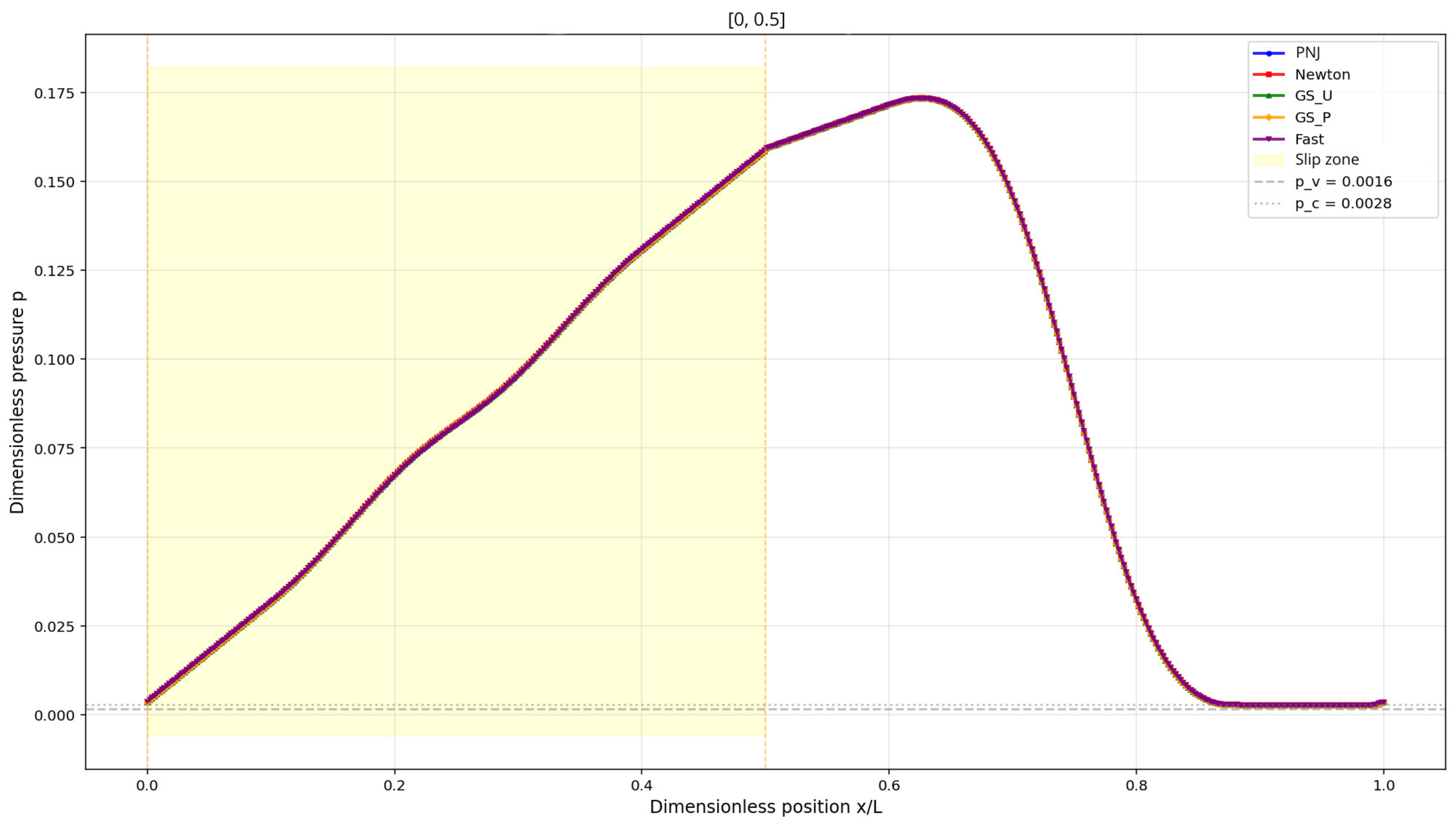}
	\hfill
	\includegraphics[width=0.48\textwidth]{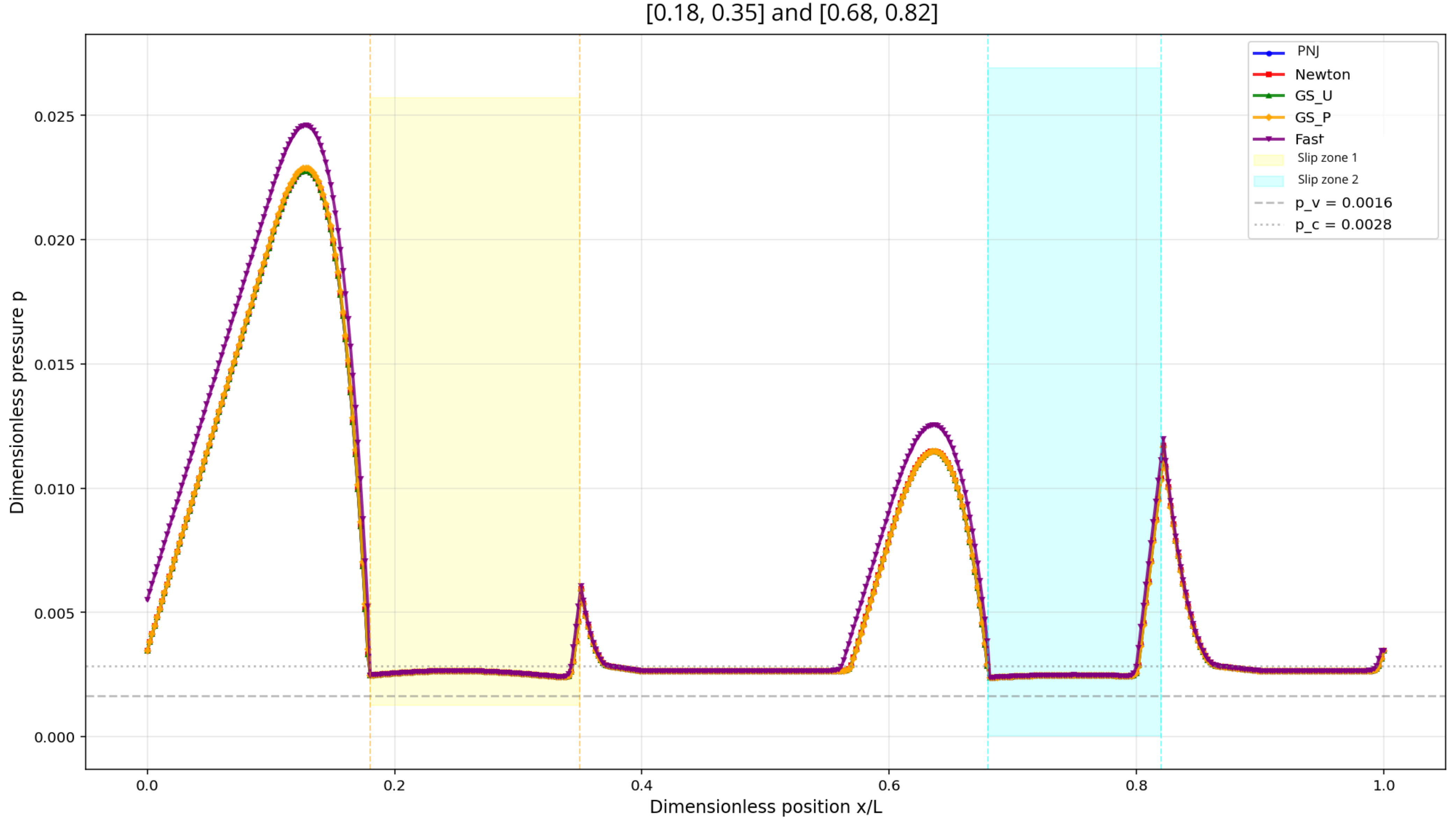}
	
	\caption{Comparison of the dimensionless pressure profiles for the four slip configurations.}
	\label{fig:pressure}
\end{figure}

\subsection{Two-dimensional problem: The twin-Parabolic Slider Bearing}
We consider the twin-parabolic slider of subsection~\ref{sec:benchmark_1D}, extended over
the physical computational domain $\Omega=[0,L_x]\times[0,L_y]$, where $L_x$ and $L_y$ denote the lengths of the lubricated surface in the $x$- and $y$-directions, respectively. To
approximate an infinitely long slider configuration, we consider a slender
geometry with $L_y=8L_x$, and the physical domain $\Omega$ is mapped onto the
dimensionless domain $\hat\Omega=[0,1]\times[0,1]$. The mesh size used throughout this study is $50 \times 50$, with a stopping tolerance of $10^{-8}$.
\subsubsection{Newton method: primal versus dual density approximation}
\label{subsec:newton_primal_dual}
We first compare Newton's method applied to the Kirchhoff-transformed problem using the primal-mesh and dual-mesh approximations of the density. The computations are performed for the slip-zone configuration $[0,0.5]\times[0,1]$.

\begin{figure}[H]
	\centering
	\begin{minipage}{0.48\textwidth}
		\centering
		\includegraphics[width=\textwidth]{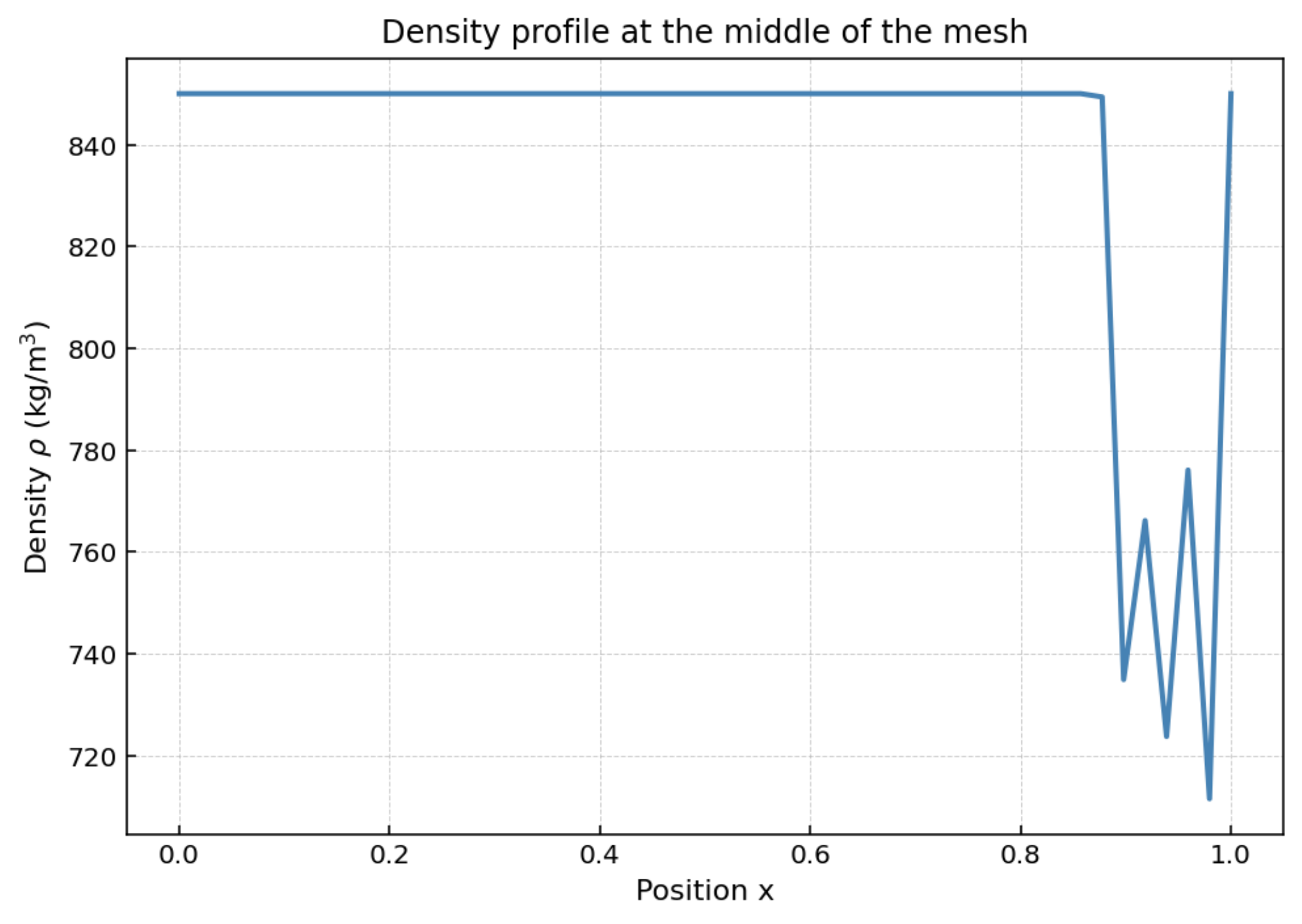}\\[2pt]
		{\small \textbf{(a)} Entire dual cell (oscillations)}
	\end{minipage}
	\hfill
	\begin{minipage}{0.48\textwidth}
		\centering
		\includegraphics[width=\textwidth]{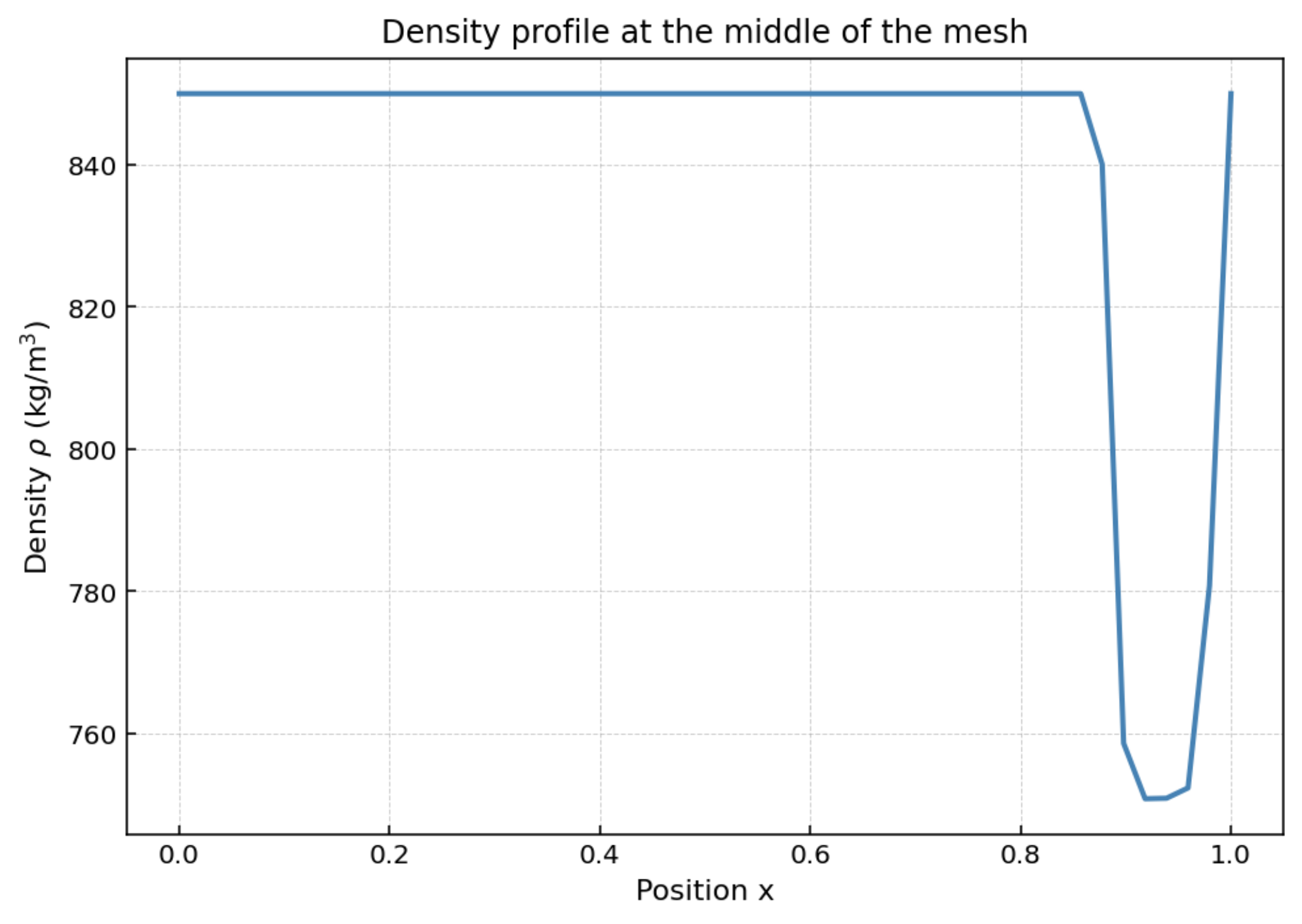}\\[2pt]
		{\small \textbf{(b)} Half dual cell (no oscillations)}
	\end{minipage}
	\caption{Comparison of the dual-mesh density approximation in the cavitation region: (a) approximation over the entire dual cell, resulting in numerical oscillations, and (b) approximation over a half dual cell, without oscillations.}
	\label{fig:density_oscillations}
\end{figure}

For the dual-mesh approximation, preliminary computations show that evaluating the density over the entire dual cell may produce oscillations in the cavitation region, as shown in Fig.~\ref{fig:density_oscillations}(a). These oscillations are not observed when a half-cell approximation is used, as illustrated in Fig.~\ref{fig:density_oscillations}(b). Therefore, the half-cell dual approximation is used in the following computations.

For the configurations considered, the convergence of Newton's method depends on the slip-zone configuration. In particular, the standard Newton method with $\omega=1$ does not converge for all the configurations tested. We also observed that a fixed relaxation parameter is not sufficient to obtain convergence for all configurations: a value of $\omega$ that works for one slip-zone configuration may lead to slow convergence or failure for another one.

To overcome this sensitivity, we employ an adaptive relaxation strategy, as considered in the literature; see, for example, Polyak and Tremba~\cite{PolyakTremba2017}. In this approach, the relaxation parameter is adjusted dynamically according to the current residual. The relaxed Newton update is written as:

\begin{equation}
	U^{(m+1)}
	=
	U^{(m)}
	+
	\omega_m\,\delta U^{(m)},
	\label{eq:newton_relax}
\end{equation}
where the relaxation parameter is defined by
\begin{equation}
	\omega_m
	=
	\min\left(
	\omega_{\max},
	\frac{C_\omega}
	{\|\mathcal{R}(U^{(m)})\|_2}
	\right).
	\label{eq:adaptive_omega}
\end{equation}
Here, $\|\mathcal{R}(U^{(m)})\|_2$ denotes the Euclidean norm of the
current residual, $C_\omega$ is a proportionality constant controlling the
step size, $\omega_{\max}$ is the maximum allowed relaxation parameter,
and $\varepsilon$ is a small regularization parameter introduced to avoid
division by zero.
Based on the convergence tests, the values
$\omega_{\max}=1$ and $C_{\omega}=6$ are used in the numerical experiments.

Thus, when the residual is large, a smaller value of $\omega_m$ is used. As the residual decreases, $\omega_m$ increases and reaches the full Newton step when $\omega_m=1$.

We then compare the primal and dual density approximations using the adaptive relaxation strategy. The results are presented in Table~\ref{tab:newton_comparison}. The primal approximation converges in $24$ iterations, while the dual approximation requires $27$ iterations. However, the CPU times differ significantly: $422.76$~s for the primal approximation compared with $1.57$~s for the dual approximation. The obtained values of the dimensionless load are $\bar{W}=0.0752$ and $\bar{W}=0.0779$, respectively, and the prescribed convergence criterion is satisfied in both cases.

For comparison, the dual approximation is also solved using the fixed relaxation parameter $\omega=0.1$. In this case, $155$ iterations are required, with a CPU time of $1.78$~s. The adaptive relaxation therefore reduces the number of Newton iterations from $155$ to $27$ for the dual approximation, while giving the same value of $\bar{W}$ up to the reported digits.
\begin{table}[H]
	\centering
	\caption{Performance comparison of Newton's method using primal and dual density approximations under fixed and adaptive relaxation strategies.}
	\label{tab:newton_comparison}
	\renewcommand{\arraystretch}{1.2}
	\begin{tabular*}{\textwidth}{@{\extracolsep\fill}lcccc}
		\toprule
		Approximation and relaxation & Iterations & CPU time (s) & $\bar{W}$ & relative error \\
		\midrule
		Primal (adaptive relaxation) & 24  & 422.76 & 0.0752 & $7.15 \times 10^{-13}$ \\
		Dual (adaptive relaxation)   & 27  & 1.57   & 0.0779 & $1.71 \times 10^{-11}$ \\
		Dual (fixed relaxation)      & 155 & 1.78   & 0.0779 & $9.20 \times 10^{-9}$  \\
		\bottomrule
	\end{tabular*}
\end{table}

\subsubsection{Computational performance}
\begin{figure}[H]
	\centering
	\includegraphics[width=0.48\textwidth]{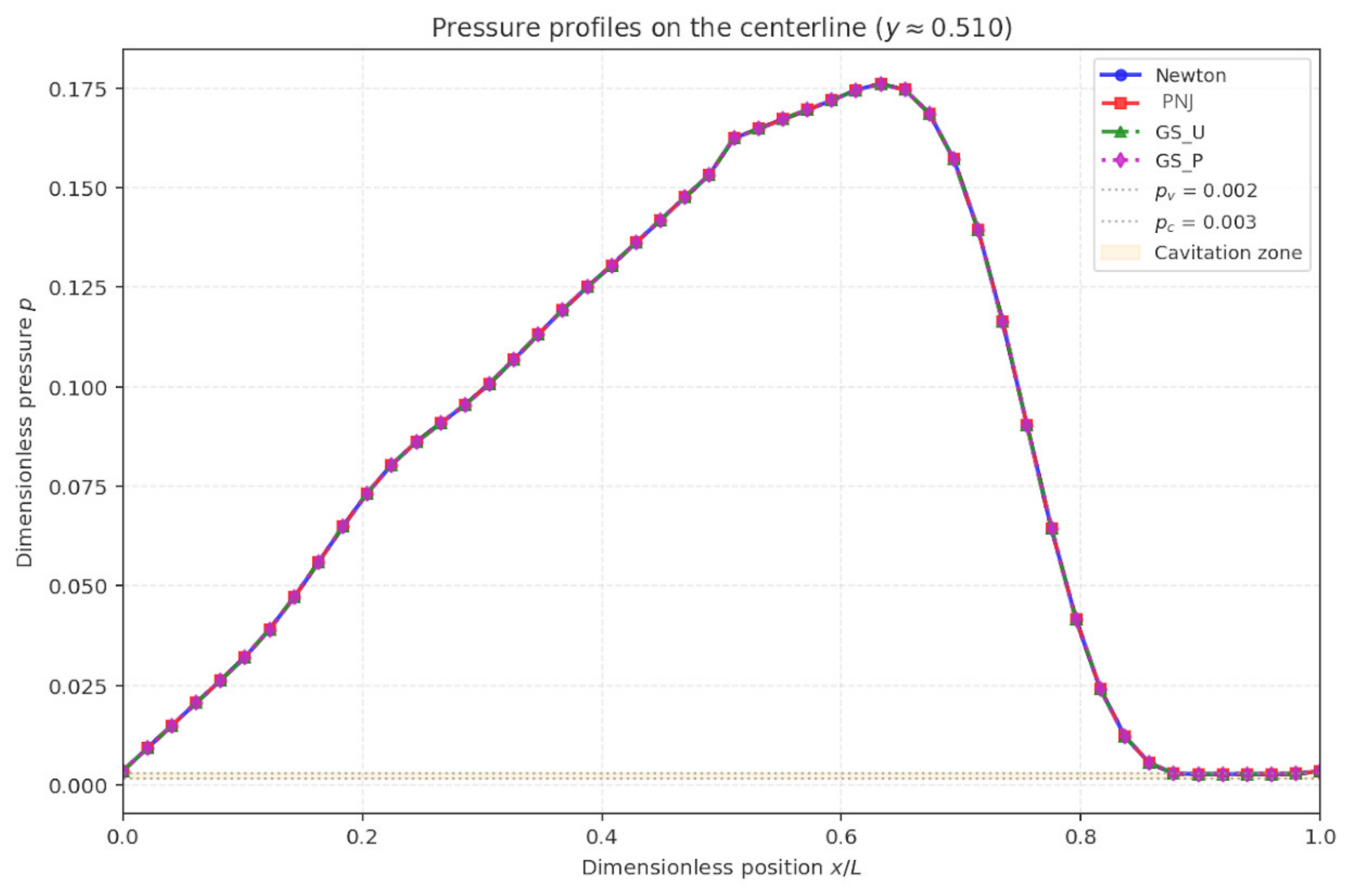}
	\caption{Pressure profile along the centerline ($y \approx 0.510$) for the slip region $[0,0.5]$.}
	\label{fig:pressure_profile_slip05}
\end{figure}

For the slip region $[0,0.5]$, the computational times show significant differences between the methods. NGSP requires the longest time ($92728.9184$~s), mainly due to its slow convergence and the large number of iterations ($5301$) required. NGSU significantly reduces the computational time to $5523.8910$~s, making it about $17$ times faster than NGSP, while requiring a comparable number of iterations ($5303$). Both Gauss--Seidel variants employed an over-relaxation parameter of $\omega = 1.83$ to accelerate convergence. Although PNJ requires a larger number of iterations ($9990$), its computational time is only $3.2777$~s, making it approximately $28{,}300$ times faster than NGSP. In comparison, Newton requires only $27$ iterations and a computational time of $1.10$~s, making it faster than PNJ method for this configuration, despite requiring significantly fewer iterations. Thus, Newton reduces the number of iterations by more than two orders of magnitude compared with Jacobi, while both methods remain computationally efficient compared with the Gauss--Seidel approaches. Despite the differences in computational cost and iteration number, all methods converge to the same solution, with $\bar{P}_{\max}=0.1761$ and $\bar{W}=0.0779$, confirming the consistency of the numerical results (Fig.~\ref{fig:pressure_profile_slip05}).

For the other slip configurations, a similar behavior is observed. Newton converges in significantly fewer iterations than PNJ, while their CPU times remain of the same order. For the no-slip configuration, PNJ and Newton require $1.10$~s and $1.15$~s, respectively, whereas for the uniform slip configuration $[0,1]$, the corresponding times are $1.06$~s and $1.20$~s. For the two separated slip regions, PNJ requires $1.18$~s compared with $1.27$~s for Newton. These results indicate that the larger number of PNJ iterations does not lead to a proportional increase in CPU time, owing to the low computational cost of each vectorized Jacobi iteration. Thus, both methods exhibit comparable computational efficiency across the different slip configurations.

\begin{table}[H]
	\centering
	\caption{Comparison of the numerical methods for different slip regions (2D case).}
	\label{tab:slip_regions}
	\renewcommand{\arraystretch}{1.2}
	\begin{tabular*}{\textwidth}{@{\extracolsep\fill}llcccc}
		\toprule
		Slip region & Method & Iterations & CPU (s) & $\bar W$ & relative error \\
		\midrule
		
		\multirow{5}{*}{No slip}
		& NGSP   & $317$   & $5695.59$ & $0.0147$  & $9.71 \times 10^{-9}$    \\
		& NGSU   & $318$   & $334.57$  & $0.0147$  & $9.59 \times 10^{-9}$     \\
		& PNJ  & $598$ & $1.10$     & $0.0147$  & $9.94 \times 10^{-9}$     \\
		& Newton  & $32$   & $1.15$   & $0.0147$ &  $2.217 \times 10^{-12}$       \\
		\midrule
		
		\multirow{5}{*}{$[0,1] \times [0,1]$}
		& NGSP   & $318$   & $5554.73$  & $0.0088$  & $9.88 \times 10^{-9}$     \\
		& NGSU   & $318$   & $324.83$  & $0.0088$   & $9.55 \times 10^{-9}$      \\
		& PNJ  & $599$    & $1.06$      & $0.0088$  & $9.927 \times 10^{-9}$      \\
		& Newton  & $48$    & $1.20$    & $0.0088$& $2.215 \times 10^{-11}$       \\
		\midrule
		
		\multirow{5}{*}{$[0,0.5]\times [0,1]$}
		& NGSP   & $5301$   & $92728.9184$ & $0.0779$  & $9.99 \times 10^{-9}$    \\
		& NGSU   & $5303$   & $5523.8910$  & $0.0779$  & $9.99 \times 10^{-9}$     \\
		& PNJ  & $9990$   & $3.28$       & $0.0779$  & $9.995 \times 10^{-9}$      \\
		& Newton  & $27$    & $1.10$      & $0.0779$  & $1.709 \times 10^{-11}$      \\
		\midrule
		
		\multirow{5}{*}{$[0.18,0.35]\times [0,1] and [0.68,0.82] \times [0,1]$}
		& NGSP   & $850 $   & $14967.48 $  & $0.0057$  & $9.85 \times 10^{-9}$     \\
		& NGSU   & $853$    & $887.27$     & $0.0057$  & $9.89 \times 10^{-9}$      \\
		& PNJ  & $1110 $  & $ 1.18$      & $0.0057$  & $9.917 \times 10^{-9}$      \\
		& Newton  & $45$    & $1.21$       & $0.0057$  & $2.41 \times 10^{-11}$        \\
		\bottomrule
	\end{tabular*}
\end{table}

\begin{figure}[H]
	\centering
	\includegraphics[width=0.48\linewidth]{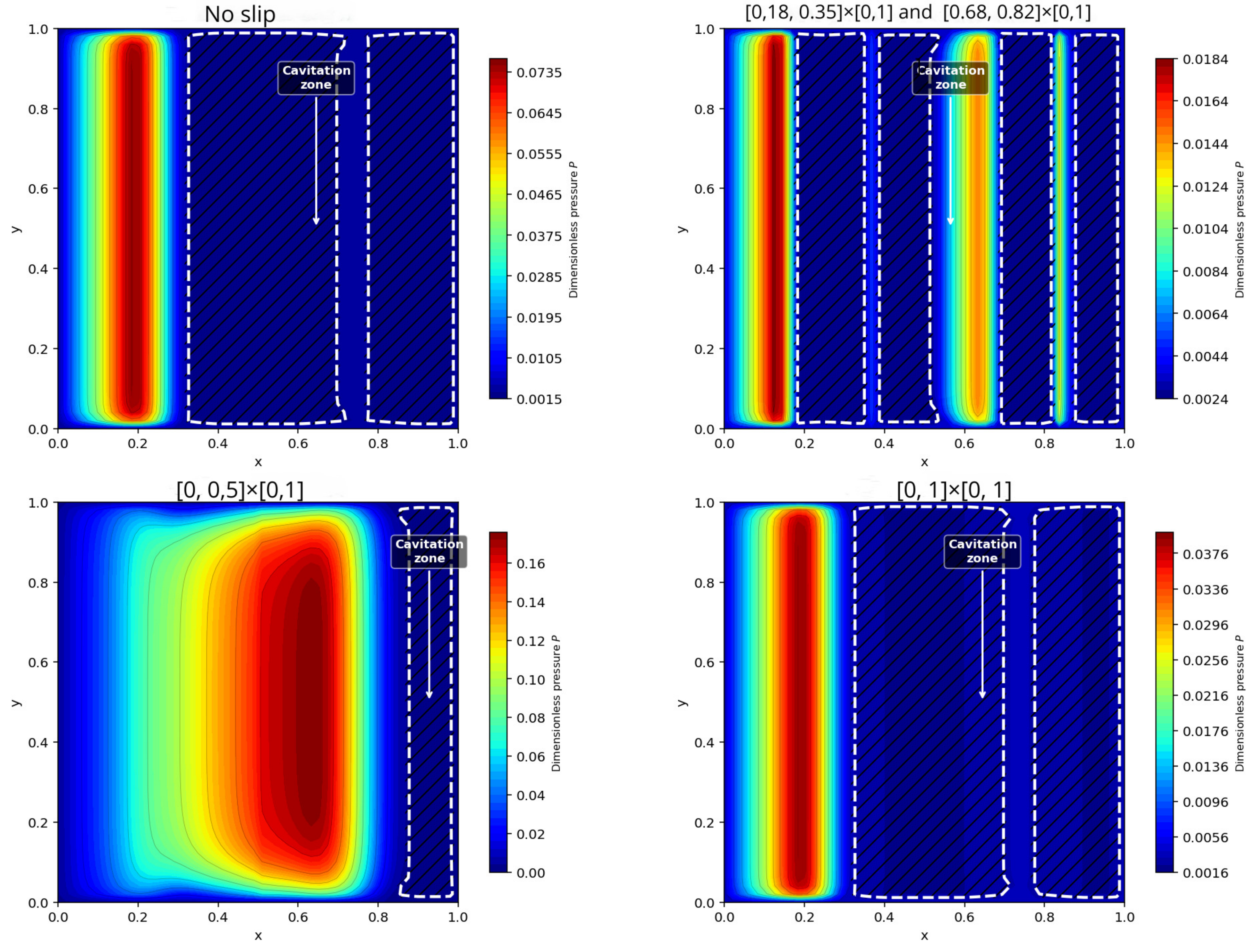}
	\caption{Pressure Distribution and Cavitation Zones in the twin parabolic slider}
	\label{ALL2B2D}
\end{figure}

\paragraph{Discussion}
\begin{itemize}
	\item Assuming a no-slip condition on the upper surface ($K(x,y)= 0$), cavitation occurs over the regions $[0.3265,0.7143]\times[0.0204,0.9796]$ and $[0.7755,0.9796]\times[0.0204,0.9796]$. The corresponding peak pressure is $0.0751$ (Fig.~\ref{ALL2B2D}, Case~1), with a load-carrying capacity of $0.0147$. As the first bump represents the stronger constriction $(h_{min1} < h_{min2})$, it governs the pressure buildup across the full width of the domain. In contrast, the shallower second bump is insufficient to generate a comparable secondary pressure peak.\\
	
	\item When slip is imposed on the upper surface over
	$c=[0.18,0.38]\times[0,1]$ and
	$[0.68,0.82]\times[0,1]$, both parabolic constrictions are partially covered by slip regions. The presence of slip modifies the local pressure distribution, leading to several pressure build-up and pressure-drop regions. Consequently, four distinct cavitation regions are observed (see Fig.~\ref{ALL2B2D}, Case~2). The pressure reaches a maximum value of $\bar P_{\max}=0.0182$, while the load-carrying capacity is $\bar W=0.0057$. Thus, this slip configuration produces weak pressure generation and results in the lowest load-carrying capacity among the four considered cases.\\
	
	\item Assuming that slip occurs on the upper surface over
	$c=[0,0.5]\times[0,1]$, the first parabolic constriction is entirely covered by the slip region, while the second constriction remains under no-slip conditions. Slip reduces the pressure build-up in the first constriction, shifting the main pressure generation to the second converging region, where the maximum pressure $\bar{P}_{\max}=0.1761$ is reached (Fig.~\ref{ALL2B2D}, Case~3). The pressure then decreases in the some regions, leading to cavitation near the outlet over
	$[0.8571,0.9796]\times[0.0204,0.9796]$. This configuration yields the highest load-carrying capacity, $\bar{W}=0.0779$, among the four cases.\\
	
	\item Assuming uniform slip over the entire upper surface, $c=[0,1]\times[0,1]$, both bumps are affected similarly, so the cavitation pattern remains close to the no-slip case. However, the pressure level is significantly reduced, with $\bar P_{\max}=0.0392$ and $\bar W=0.0088$ (Fig.~\ref{ALL2B2D}, Case~4). Thus, uniform slip weakens the wedge effect throughout the domain without significantly altering the overall pressure distribution.
\end{itemize}

To investigate the influence of the bearing dimensions, we consider four configurations with \(L_y/L_x=1,2,4,\) and \(8\), using two separated slip regions, \([0.18,0.35]\times[0,1]\) and \([0.68,0.82]\times[0,1]\). As shown in Fig.~\ref{fig:pressure}, the pressure distributions are qualitatively similar for all configurations. In particular, the case \(L_y/L_x=8\) approximates an infinitely long bearing, whereas \(L_y/L_x=1\) corresponds to a square bearing with \(L_x=L_y\). The dimensionless load decreases slightly from \(0.0058\) for the long-bearing configuration to \(0.0054\) for the square configuration, representing a reduction of approximately \(6.9\%\). In terms of computational performance, Newton's method requires 33 iterations and 1.31 s for \(L_y/L_x=8\), compared with 37 iterations and 1.34 s for \(L_y/L_x=1\). Thus, the square configuration yields a slightly lower load-carrying capacity, while the CPU times remain nearly identical.

\begin{table}[H]
	\centering
	\caption{Performance of Newton's method for different bearing
		aspect ratios with two separated slip regions.}
	\label{tab:aspect_ratio}
	\renewcommand{\arraystretch}{1.2}
	\begin{tabular*}{\textwidth}{@{\extracolsep\fill}ccccc}
		\toprule
		$L_y/L_x$ & Iterations & CPU (s) & $\bar W$ & relative error \\
		\midrule
		8 & 33 & 1.31 & 0.0058 & $1.568\times10^{-11}$ \\
		4 & 34 & 1.47 & 0.0057 & $2.982\times10^{-14}$ \\
		2 & 35 & 1.50 & 0.0056 & $4.519\times10^{-10}$ \\
		1 & 37 & 1.34 & 0.0054 & $1.704\times10^{-10}$ \\
		\bottomrule
	\end{tabular*}
\end{table}

\begin{figure}[H]
	\centering
	
	\includegraphics[width=0.48\textwidth]{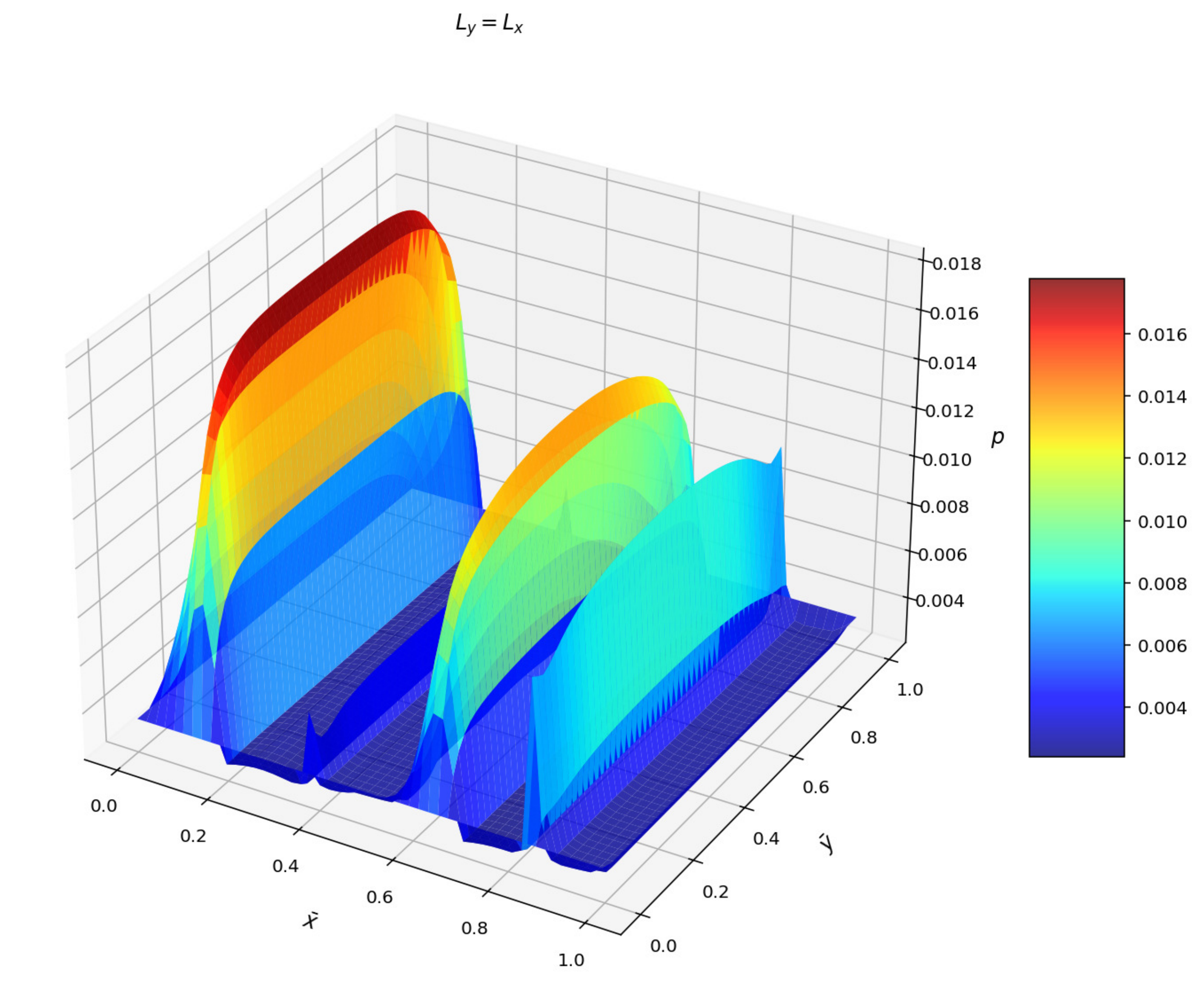}
	\hfill
	\includegraphics[width=0.48\textwidth]{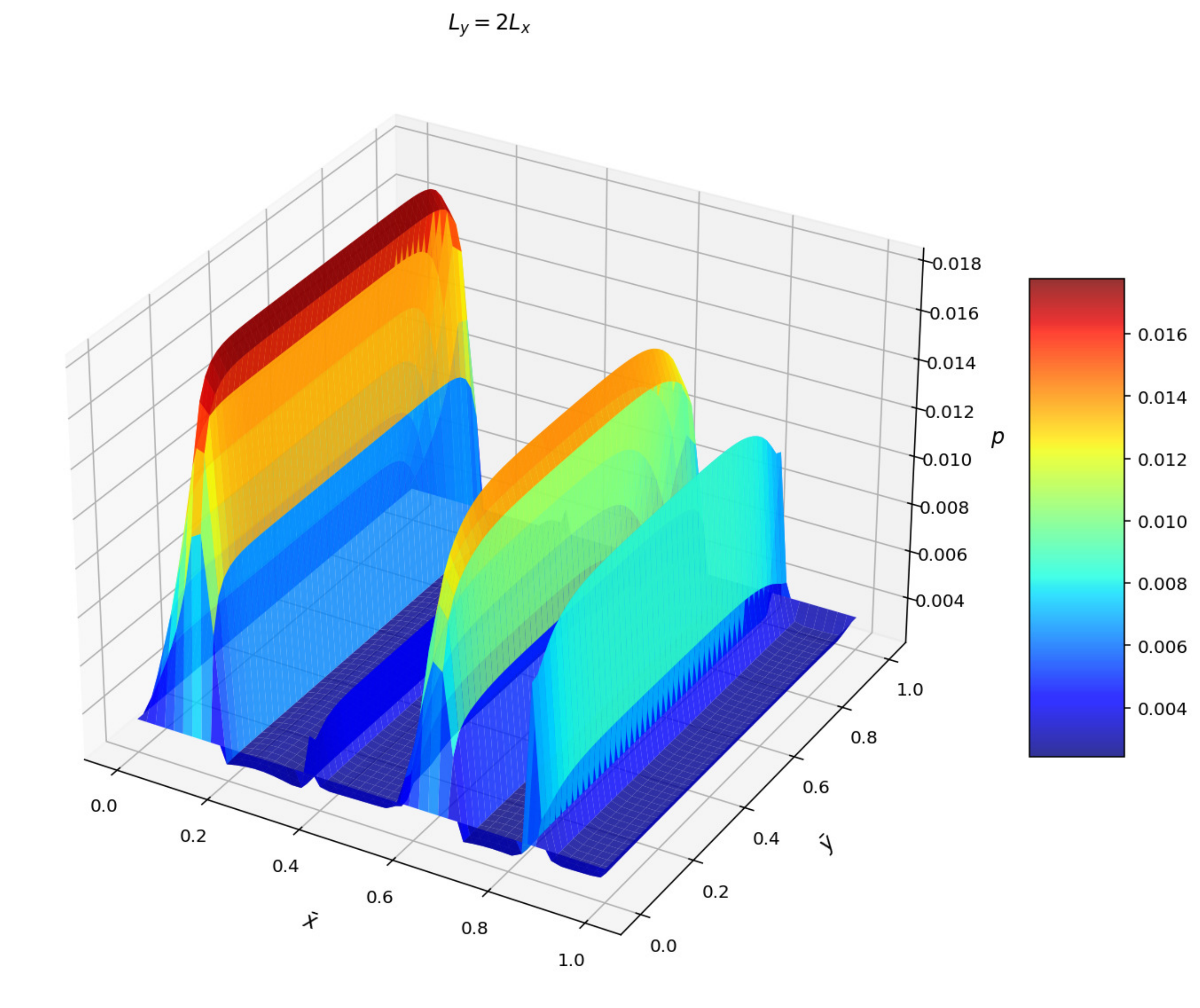}
	
	\vspace{0.4cm}
	
	\includegraphics[width=0.48\textwidth]{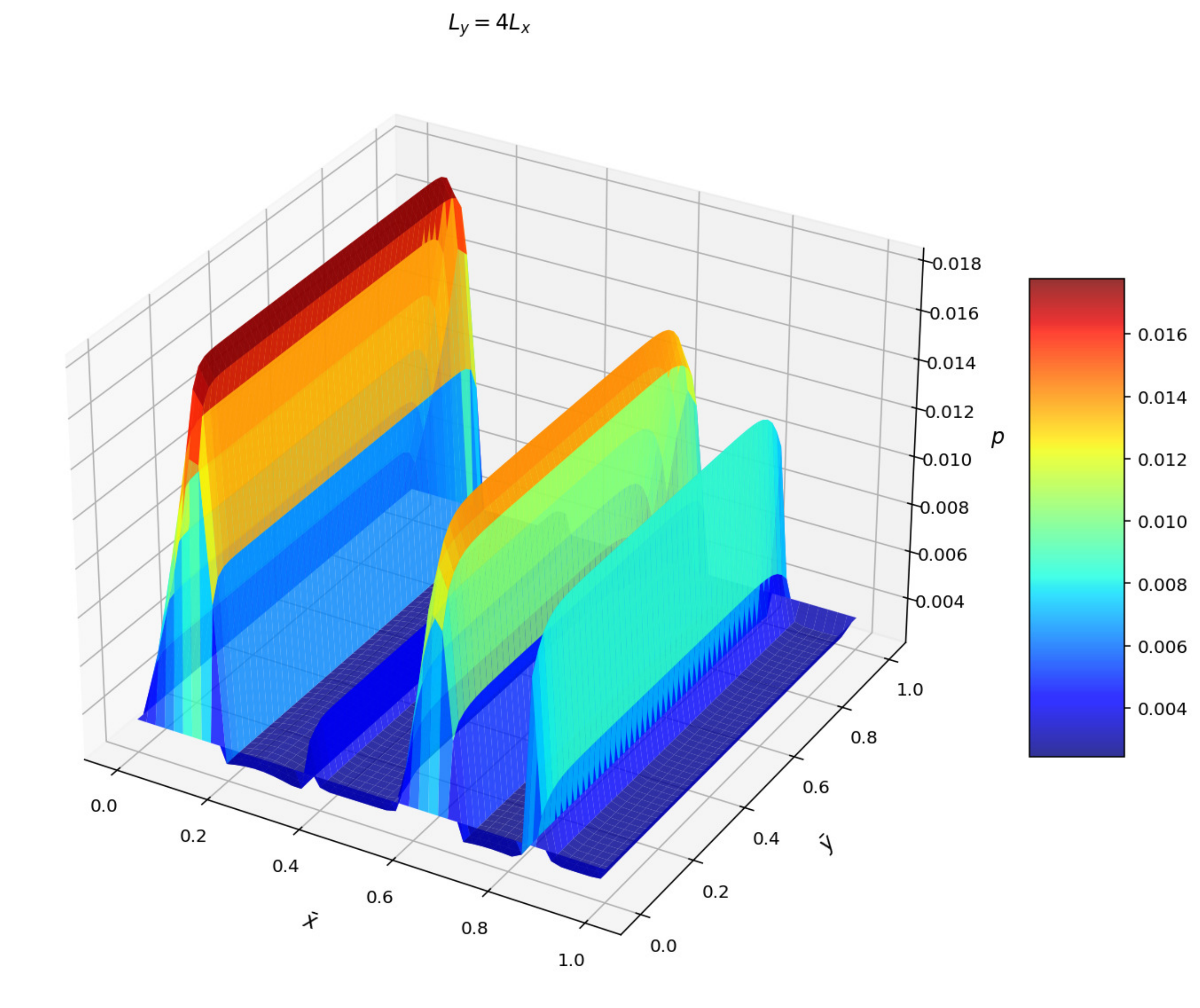}
	\hfill
	\includegraphics[width=0.48\textwidth]{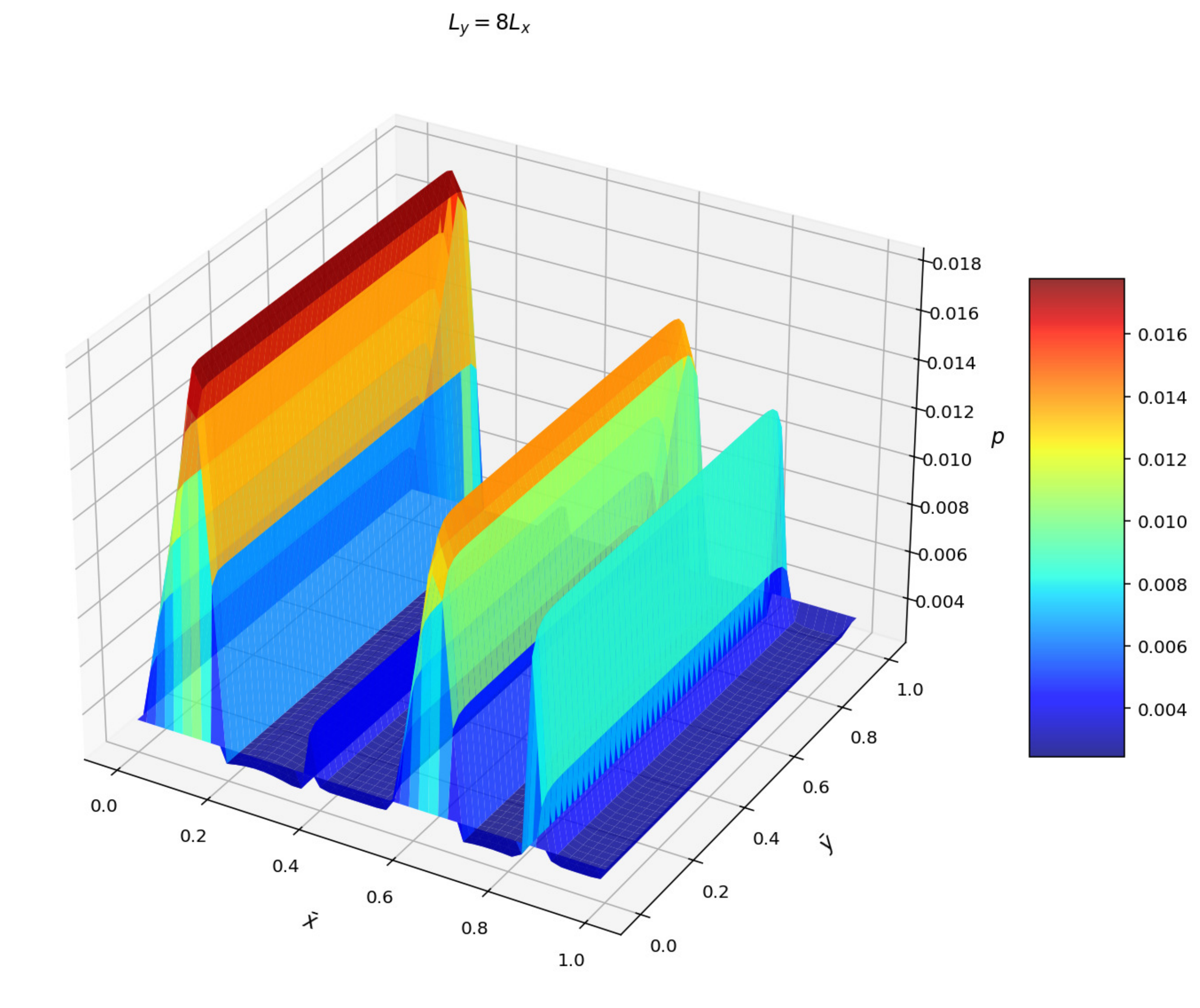}
	
	\caption{Dimensionless pressure distributions for different
		bearing aspect ratios with heterogeneous slip regions.}
	\label{fig:pressure}
\end{figure}

\subsection{Two-dimensional problem: The textured surface}
\label{subsec:slip_texture}
\begin{figure}[H]
	\centering
	\begin{minipage}{0.48\textwidth}
		\centering
		\includegraphics[width=\textwidth]{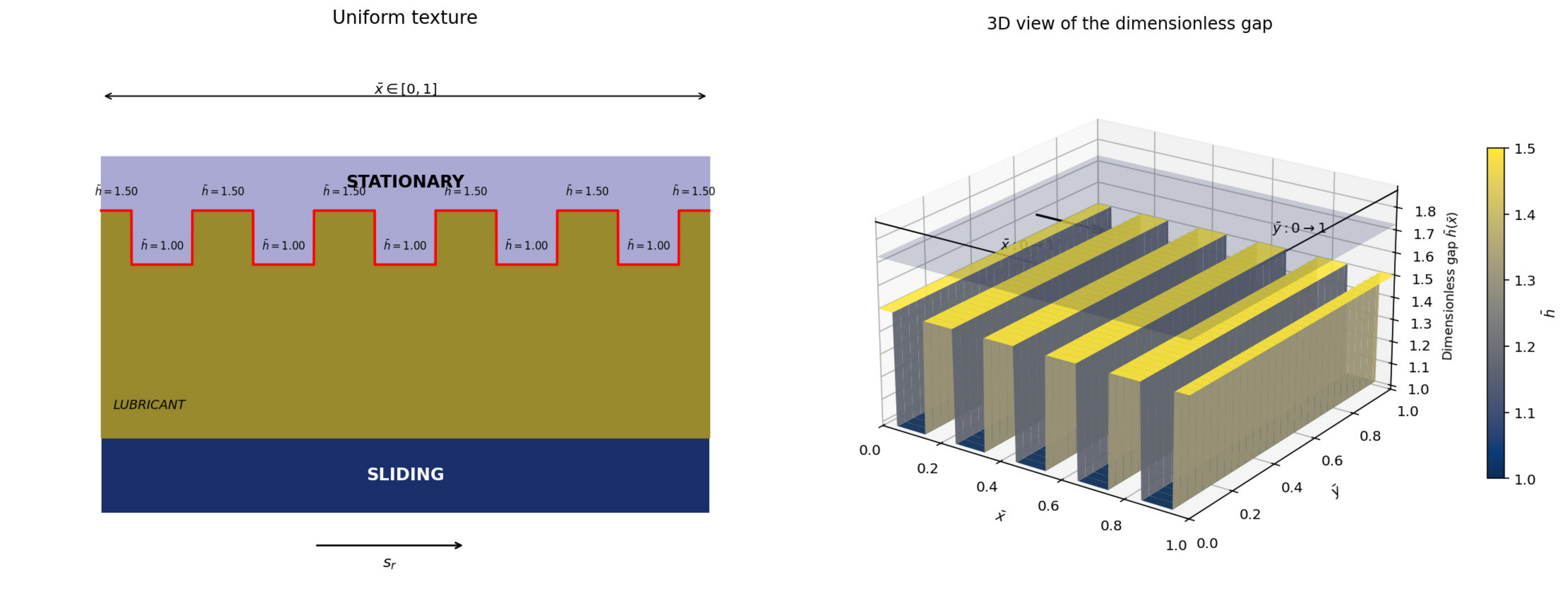}\\[2pt]
		{\small \textbf{(a)} Uniform texture}
	\end{minipage}
	\hfill
	\begin{minipage}{0.48\textwidth}
		\centering
		\includegraphics[width=\textwidth]{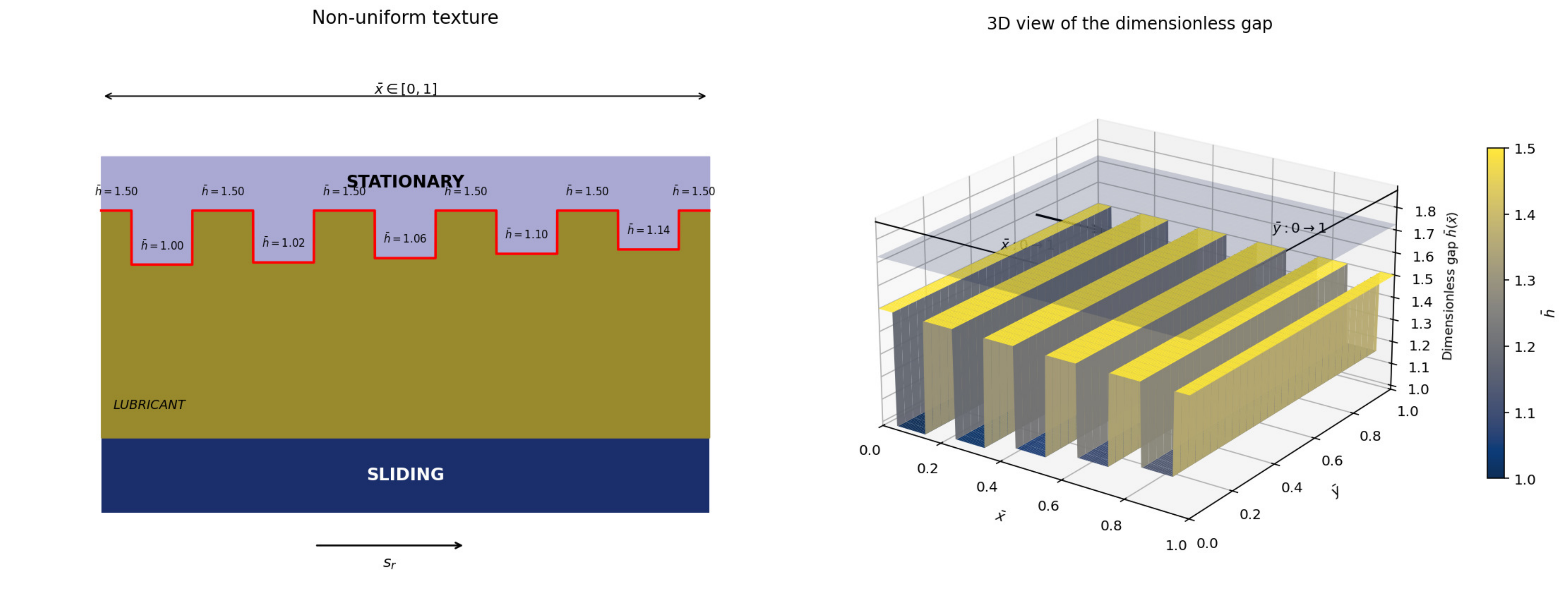}\\[2pt]
		{\small \textbf{(b)} Non-uniform texture with progressively increasing recess depth.}
	\end{minipage}
	\caption{Dimensionless textured film-thickness profile $\bar{H}(\bar{X})$.}
	\label{fig:textured_gap}
\end{figure}

We next investigate the effect of boundary slip on a textured lubricated
contact. The textured geometry considered here is taken from the gap
introduced by Ausas~\citep{Ausas2007}. In dimensionless form, the film
thickness is a piecewise-constant function $\bar{H}(\bar{X})$,
alternating between $\bar{H} = 1$ and $\bar{H} = 1.5$, as illustrated in
Fig.~\ref{fig:textured_gap}(a). We also consider a second, irregular texture, in which the film
thickness in the thinner cells increases progressively along the
$\bar{x}$-direction, taking the values $\bar{H} = 1.00,\,1.02,\,1.06,\,1.10,\,1.14$,
while the thicker cells are kept fixed at $\bar{H} = 1.5$, as shown
in Fig.~\ref{fig:textured_gap}(b).

The parameters used for the numerical study of the textured
lubricated contact are summarized in Table~\ref{T3}.
\begin{table}[H]
	\centering
	\caption{Physical data.}
	\label{T3}
	\renewcommand{\arraystretch}{1.2}
	\begin{tabular}{l l}
		\toprule
		\textbf{Parameter} & \textbf{Value} \\
		\midrule
		Sliding velocity $sr$ (m/s) & 0.3 \\
		Viscosity $\mu$ (Pa s) & 0.04 \\
		Length $L_x$ (m) & 1 \\
		Length $L_y$ (m) & $8L_x$ \\
		$h_0$  & 1.5 \\
		Slip coefficient $\alpha$ & 1 \\
		Boundary fluid pressure $\bar{p}_{ext}$  & $1 \times 10^{-3}$ \\
		\bottomrule
	\end{tabular}
\end{table}

The pressure-density relationship is governed by the following parameters: 
\[
\rho_l = 850 \, \text{kg/m}^3, \quad \rho_v = 0.017 \, \text{kg/m}^3, \quad p_c = 8.2 \times 10^{4} \, \text{Pa}, \quad p_v = 4.7 \times 10^{4} \, \text{Pa}
\]

\subsubsection{influence of the cavitation model}
In \citep{Ausas2007}, the influence of the cavitation model on the performance of lubricated devices was demonstrated. It is therefore convenient to compare the fully compressible cavitation model of \citep{Bayada2013} with the classical Swift--Stieber model \citep{Christopherson1941}. Although the latter does not conserve mass, it is widely used because of its numerical simplicity. In this approach, the density is assumed to be constant, and cavitation is described by the pressure constraint $p\geq p_c$. When the pressure reaches the cavitation threshold $p_c$, the pressure is maintained at this value in the cavitated region. 

\begin{figure}[H]
	\centering
	\begin{minipage}{0.48\textwidth}
		\centering
		\includegraphics[width=\textwidth]{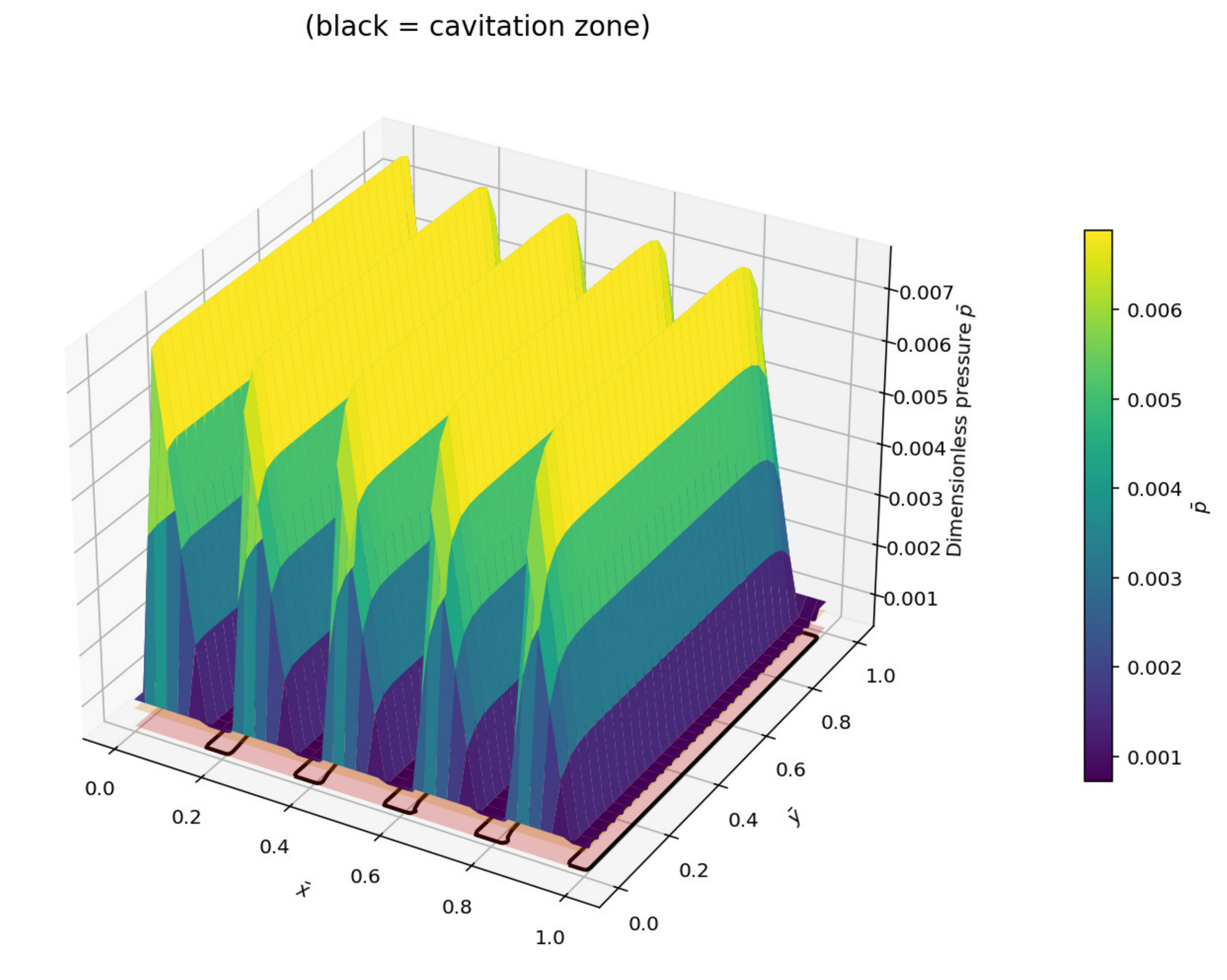}\\[2pt]
		{\small \textbf{(c)} Dimensionless pressure for uniform texture depth (3D).}
	\end{minipage}
	\hfill
	\begin{minipage}{0.48\textwidth}
		\centering
		\includegraphics[width=\textwidth]{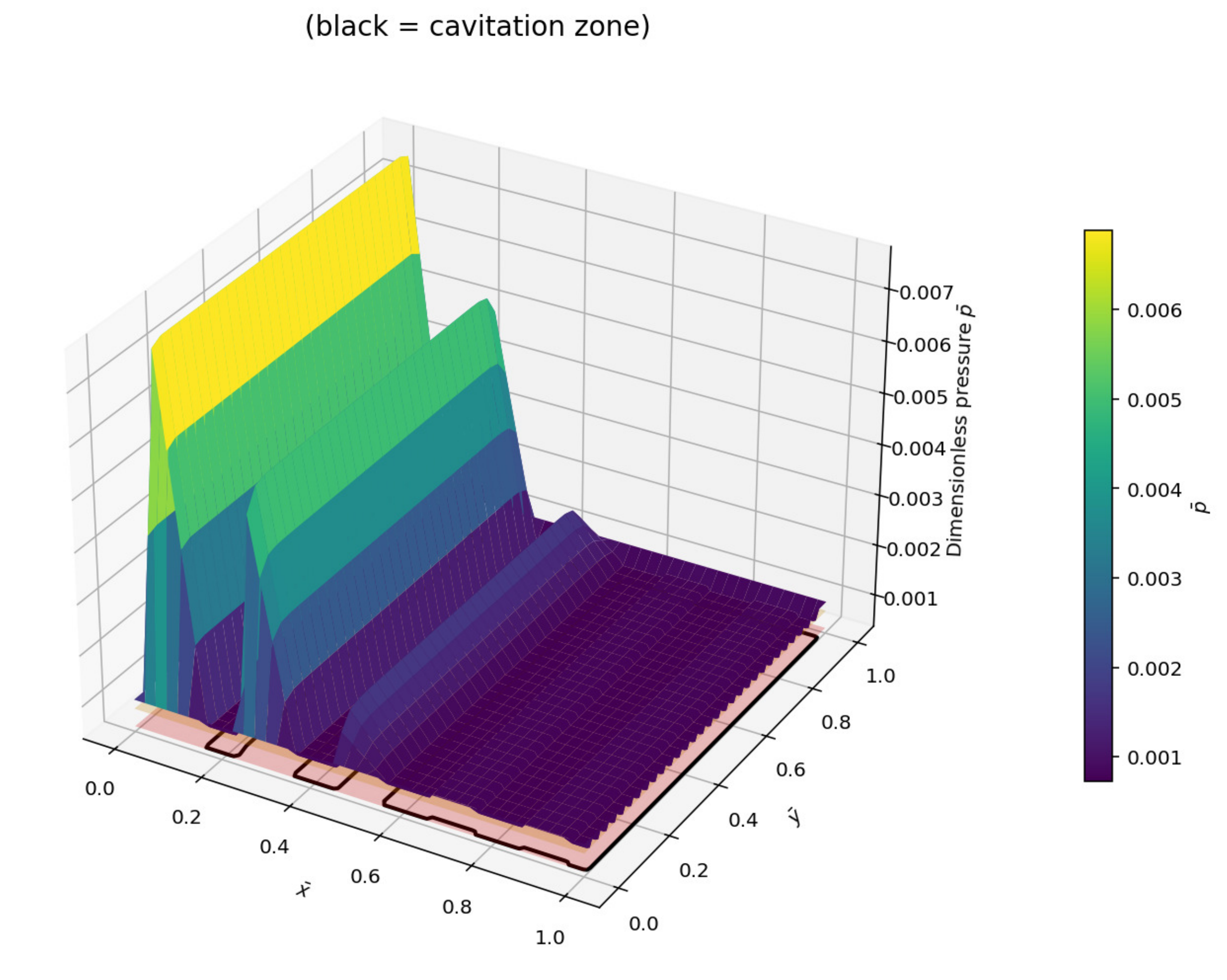}\\[2pt]
		{\small \textbf{(d)} Dimensionless pressure for non-uniform texture depth (3D).}
	\end{minipage}
	\caption{Dimensionless pressure distributions for the fully compressible model.}
	\label{fig:textured_RC}
\end{figure}

\begin{figure}[H]
	\centering
	\begin{minipage}{0.48\textwidth}
		\centering
		\includegraphics[width=\textwidth]{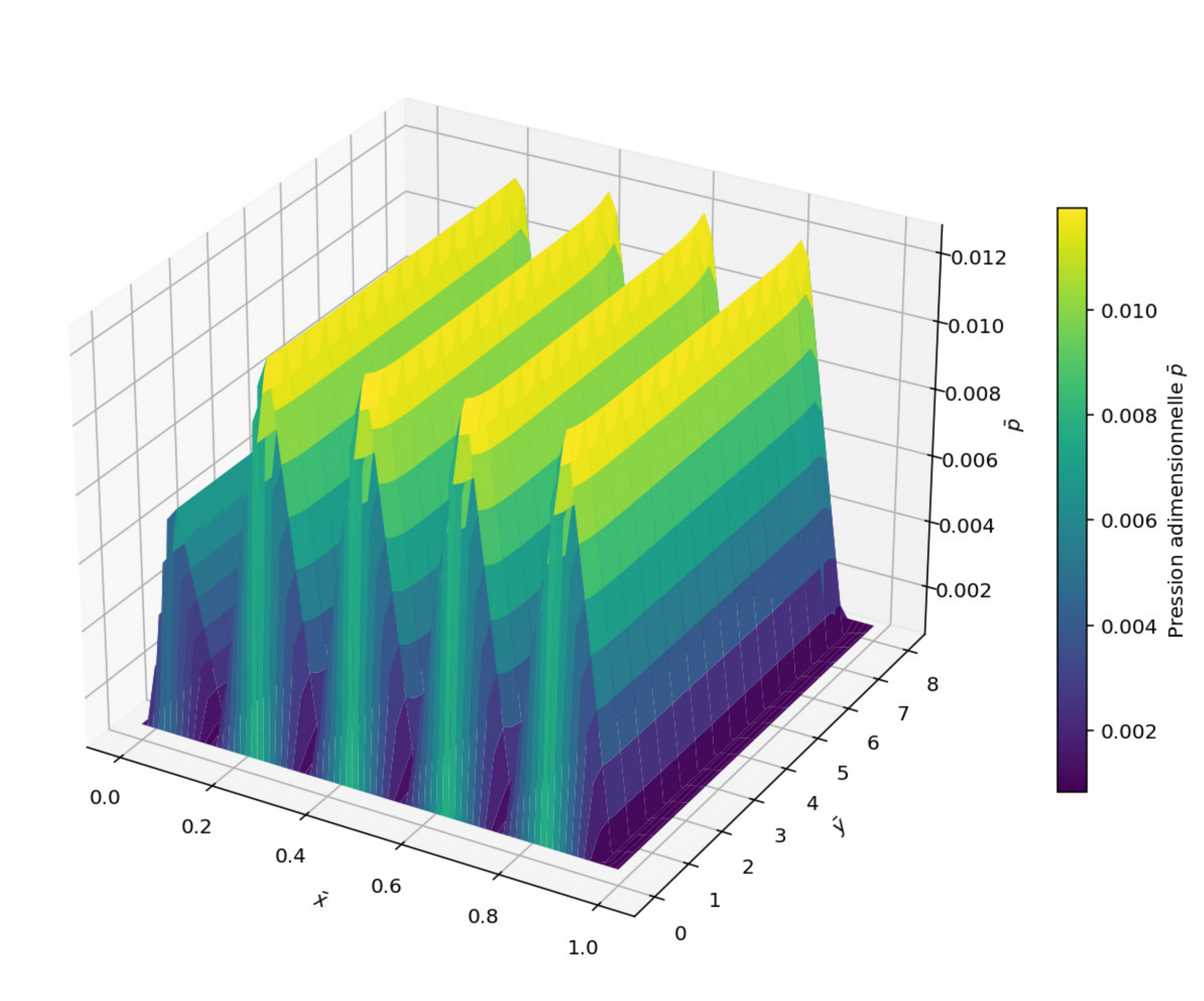}\\[2pt]
		{\small \textbf{(e)} Dimensionless pressure for uniform texture depth (3D).}
	\end{minipage}
	\hfill
	\begin{minipage}{0.48\textwidth}
		\centering
		\includegraphics[width=\textwidth]{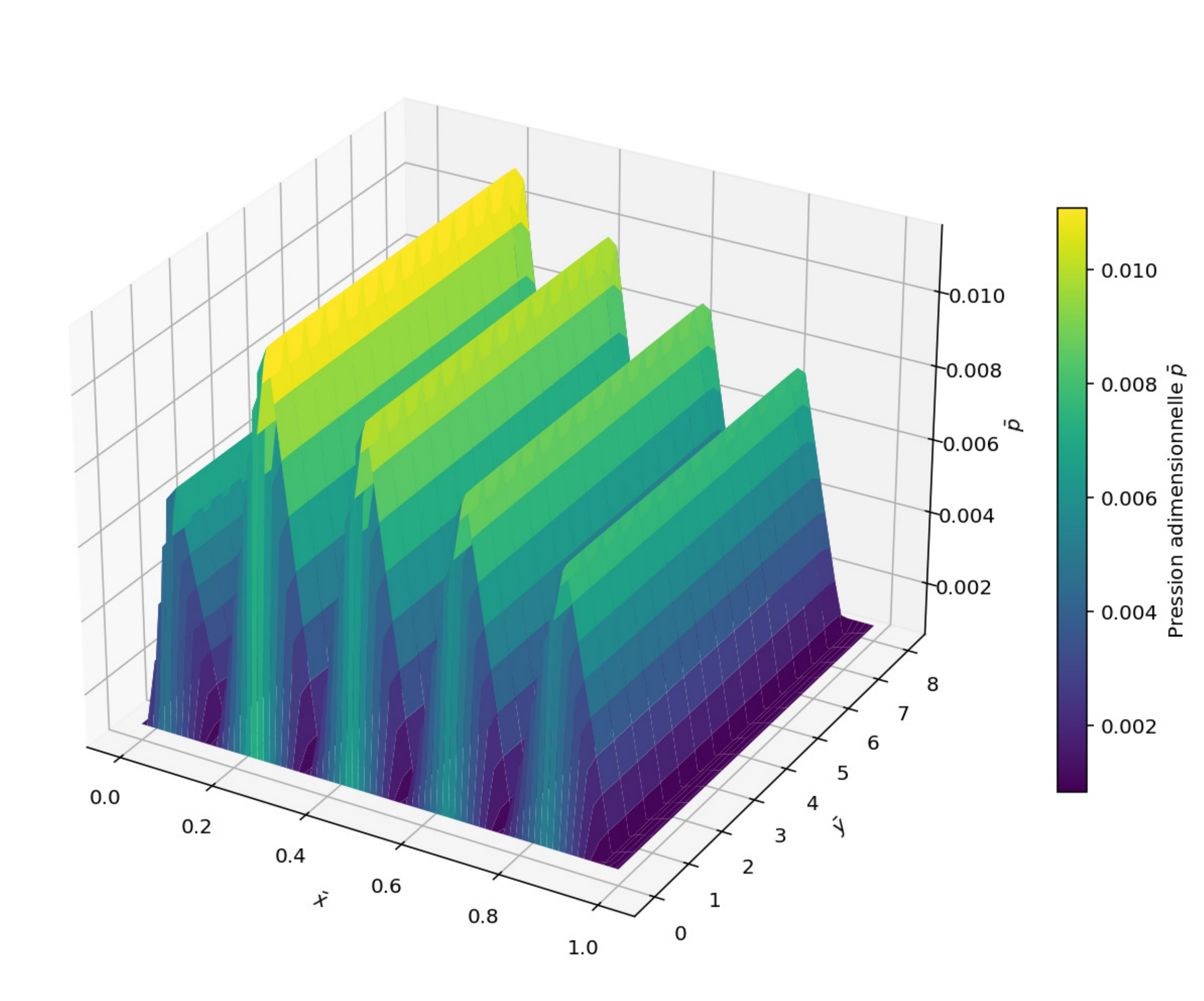}\\[2pt]
		{\small \textbf{(f)} Dimensionless pressure for non-uniform texture depth (3D).
		}
	\end{minipage}
	\caption{Dimensionless pressure distributions for Swift Steiber Model.}
	\label{fig:textured_SS}
\end{figure}

The comparison is performed for the textured configurations shown in Fig.~\ref{fig:textured_gap}(a) and (b). Since Ausas et al.~\citep{Ausas2007} considered a one-dimensional configuration, we set $L_y=8L_x$ to approximate an infinitely long lubricated device. The computations are carried out using Newton's method on a $160\times 32$ mesh, with a convergence tolerance of $10^{-8}$. The corresponding pressure distributions obtained with the fully compressible model are shown in
Fig.~\ref{fig:textured_RC} (c) and (d), while those obtained with the
Swift--Stieber model are presented in Fig.~\ref{fig:textured_SS} (e) and (f).

The Swift--Stieber model predicts significantly higher load capacities than
the fully compressible model, with $\bar{W}=0.0058$ and $0.0048$ for the uniform
and non-uniform textures, respectively, compared with $\bar{W}=0.0032$ and
$0.0016$ for the fully compressible model, as reported in
Table~\ref{tab:load_RC_SS}. The Swift--Stieber model is a non-conservative
cavitation model, whereas the fully compressible model is mass-conserving
and accounts for the liquid--mixture--vapor transition. For the Swift–Stieber model, the pressure distributions in Fig.~\ref{fig:textured_SS} show pressure
recovery at the texture cells, whereas the fully compressible model
predicts more extended cavitated regions, particularly for the
non-uniform texture, as shown in Fig.~\ref{fig:textured_RC}. The results obtained are consistent with the findings of Ausas et al.~\citep{Ausas2007}, who reported similar discrepancies between Swift--Stieber
model and the mass-conserving Elrod--Adams (JFO) model.

\begin{table}[H]
	\centering
	\caption{Load capacity $\bar{W}$ for the fully compressible and Swift--Stieber cavitation models, in the uniform and non-uniform texture cases.}
	\label{tab:load_RC_SS}
	\renewcommand{\arraystretch}{1.2}
	\begin{tabular}{lcc}
		\toprule
		\multirow{2}{*}{Model} & \multicolumn{2}{c}{Load capacity $W$} \\
		\cmidrule(lr){2-3}
		& Uniform texture & Non-uniform texture \\
		\midrule
		The fully compressible model   & $0.0032$ & $0.0016$\\
		Swift--Stieber                 & $0.0058$ &$ 0.0048 $ \\
		\bottomrule
	\end{tabular}
\end{table}

\subsubsection{influence of the location of slip zones}
To assess the effect of the slip distribution, five configurations are considered on the upper surface. The first corresponds to the no-slip case. The second configuration corresponds to a fully textured surface with slip, $[0,1]\times[0,1]$. The third configuration considers slip over the upstream half, $[0,0.5]\times[0,1]$. The fourth configuration considers slip over the central region, $[0.25,0.75]\times[0,1]$, while the fifth configuration includes two separated slip regions near the inlet and outlet, $[0,0.25]\times[0,1]$ and $[0.75,1]\times[0,1]$.

The resulting dimensionless maximum pressure $\bar{P}_{\max}$, load-carrying capacity $\bar{W}$, CPU time, number of Newton iterations, and final error are reported in Table~\ref{tab:slip_texture}. The corresponding dimensionless pressure distributions are shown in Fig.~\ref{RCUNIFT}.

The effect of slip on the load-carrying capacity depends strongly on its spatial distribution. Full slip, $c=[0,1]\times[0,1]$, slightly reduces the load from $\bar W=0.0032$ to $\bar W=0.0031$ for the textured surface, whereas for other surface geometries, such as the twin-parabolic bearing, the same full-slip configuration can reduce the load-carrying capacity to approximately half of the corresponding no-slip value (see Table~\ref{tab:slip_regions}). Slip over $[0,0.5]\times[0,1]$ increases it to $\bar W=0.0040$. Slip over the central region, $[0.25,0.75]\times[0,1]$, further reduces the load to $\bar W=0.0025$, whereas the two separated slip regions, $[0,0.25]\times[0,1]$ and $[0.75,1]\times[0,1]$, give the highest load, $\bar W=0.0041$. Thus, among the configurations considered, placing slip in two separated regions near the inlet and outlet provides the highest load, whereas central slip reduces the load.

From a numerical point of view, Newton's method converges for all configurations. The number of iterations ranges from 46 to 74, while the CPU time varies between 5.23 and 5.48~s. The final errors remain below $10^{-8}$ for all configurations, confirming the convergence of the Newton solver for the prescribed slip distributions.

\begin{table}[H]
	\centering
	\caption{Effect of the Navier slip distribution on the textured surface.}
	\label{tab:slip_texture}
	\renewcommand{\arraystretch}{1.2}
	\begin{tabular}{l c c c c c}
		\toprule
		\textbf{Slip region} $c$ & $\bar P_{\max}$ & \text{Load} $\bar W$ & \text{CPU time} (s) & \text{Iterations} & \text{Error} \\
		\midrule
		No slip                & $0.0078$ & $0.0032$ & 5.23 & 49 & $1.035\times10^{-11}$ \\
		$[0, 1] \times [0, 1]$  & $0.0074$ & $0.0031$ & 5.31 & 51 & $1.5\times10^{-11}$ \\
		$[0,0.5]\times [0,1]$         & $0.0105 $ & $0.0040$ & 5.35 & 52 & $1.28\times10^{-12}$ \\
		$c=[0.25,0.75] \times [0,1]$,           & $0.0078$ & $0.0025$ & 5.48 & 74 & $1.48\times10^{-14}$ \\
		$c=[0,0.25] \times [0,1]$ and $[0.75,1] \times [0,1]$  & $0.0105$ & $0.0041$ & 5.27 & 46 & $6.25\times10^{-9}$ \\
		\bottomrule
	\end{tabular}
\end{table}

\begin{figure}[H]
	\centering
	\includegraphics[width=0.55\textwidth]{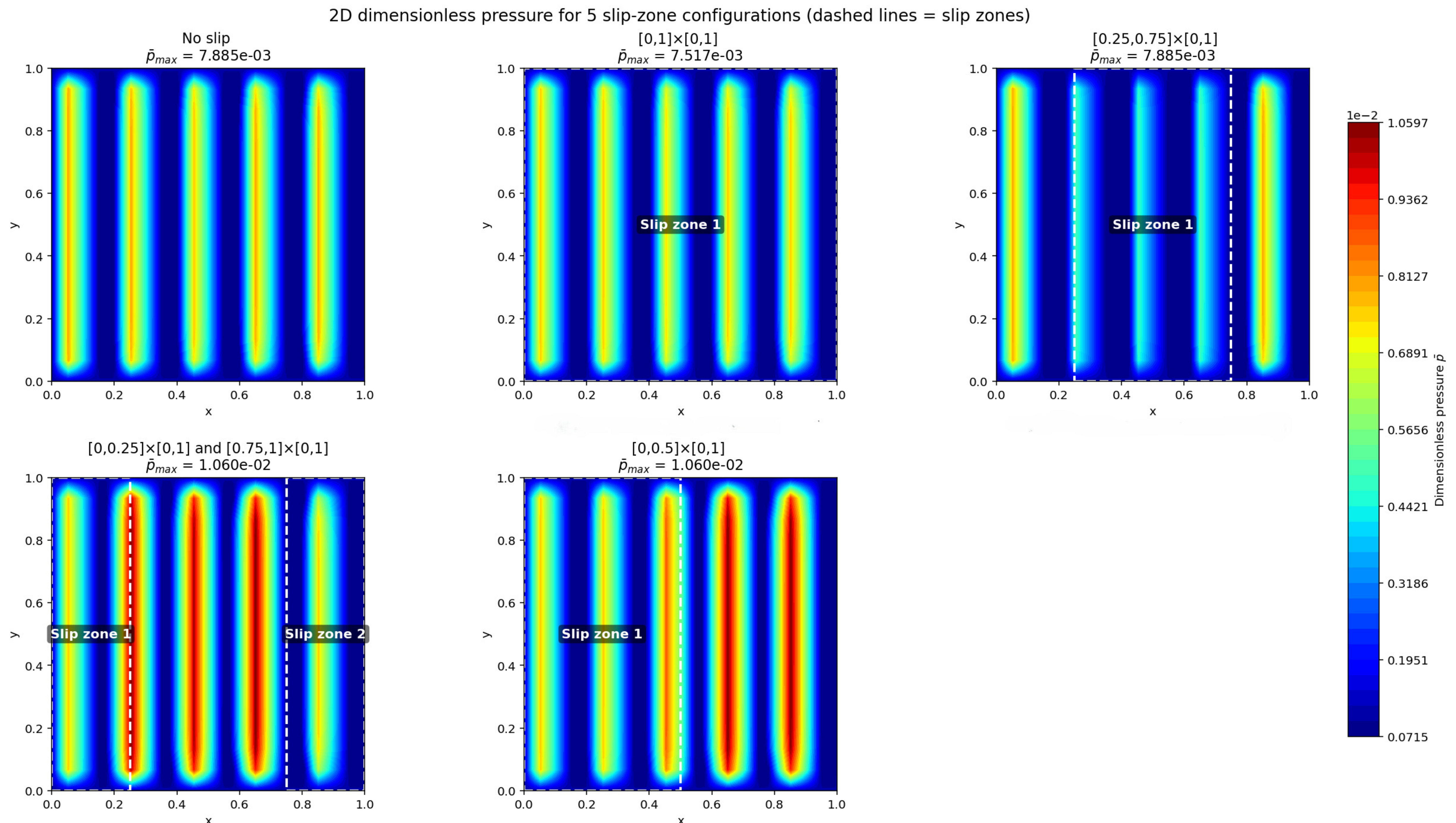}
	\caption{Dimensionless pressure distributions for five slip configurations with uniform texture depth.}
	\label{RCUNIFT}
\end{figure}

\section{Conclusion}
This work presented an efficient numerical framework for solving the one and two-dimensional compressible Reynolds equation with Navier slip boundary and mass-conserving cavitation. By using the Kirchhoff transformation and a dual-mesh finite element discretization, we successfully decoupled the nonlinear diffusion term, enabling robust and fast iterative solvers. Our comparative analysis demonstrates that while the proposed parallelized nonlinear Jacobi method is highly efficient and significantly faster than classical Gauss–Seidel approaches, Newton’s method with adaptive relaxation remains the fastest solver overall. This highlights the superior convergence rate of Newton's method, although the vectorized Jacobi scheme offers a compelling alternative due to its simplicity and parallelization potential. Furthermore, numerical benchmarks confirm that heterogeneous Navier slip patterns can significantly enhance load-carrying capacity when strategically positioned, unlike uniform slip which tends to degrade performance. 

\appendix
\section{Algorithms}
\label{app:algorithms}

\newcounter{algo}

\begin{figure}[H]
	\renewcommand{\figurename}{Algorithm}
	\stepcounter{algo}
	\renewcommand{\thefigure}{\arabic{algo}}
	\renewcommand{\theHfigure}{alg.\arabic{algo}}
	\small
	\begin{enumerate}
		\item[] \textbf{Input:} $U^{(0)}$, matrices $M_{\mathrm{off}}$, $E_{\mathrm{off}}$, tolerance $\varepsilon$
		\item $m \gets 0$
		\item \textbf{repeat}
		\begin{enumerate}
			\item Compute $T_{1_i}^{(m-1)}$ for all $i$ \hfill (Eq.~\eqref{eq:T1_vectorized})
			\item Compute $U_{1_i}$ and $U_{2_i}$ for all $i$ \hfill (Eqs.~\eqref{eq:prop1}--\eqref{eq:prop2})
			\item Compute $\alpha_i, \beta_i, \gamma_i$ for all $i$ \hfill (Eqs.~\eqref{eq:alpha_i}--\eqref{eq:gamma_i})
			\item Compute $U_{3_i}$ for all $i$ \hfill (Eq.~\eqref{eq:prop3})
			\item \textbf{for} each node $i$ (vectorized selection, without an explicit loop) \textbf{do}
			\begin{enumerate}
				\item \textbf{if} $U_{1_i}+R(\bar p_{\mathrm{ext}}) \le \bar u_v$ \textbf{then} $U_i^{(m)} \gets U_{1_i}$
				\item \textbf{else if} $U_{2_i}+R(\bar p_{\mathrm{ext}}) \ge \bar u_c$ \textbf{then} $U_i^{(m)} \gets U_{2_i}$
				\item \textbf{else} $U_i^{(m)} \gets U_{3_i}$
			\end{enumerate}
			\item \textbf{end for}
			\item Compute the relative error $e^{(m)}=\|U^{(m)}-U^{(m-1)}\|/\|U^{(m)}\|$
			\item $m \gets m+1$
		\end{enumerate}
		\item \textbf{until} $e^{(m)}<\varepsilon$
		\item \textbf{return} $U^{(m)}$
	\end{enumerate}
	\caption{Parallelized nonlinear Jacobi algorithm (PNJ)}
	\label{alg:jacobi_vectorized}
	\addtocounter{figure}{-1}
\end{figure}

\begin{figure}[H]
	\renewcommand{\figurename}{Algorithm}
	\stepcounter{algo}
	\renewcommand{\thefigure}{\arabic{algo}}
	\renewcommand{\theHfigure}{alg.\arabic{algo}}
	\small
	\begin{enumerate}
		\item[] \textbf{Input:} $U^{(0)}$, matrices $M$, $E$, tolerance $\varepsilon$, relaxation parameter $\omega$
		\item $m \gets 0$
		\item \textbf{repeat}
		\begin{enumerate}
			\item \textbf{for} $i = 1, \ldots, n_v$ \textbf{do} \hfill (sequential forward sweep)
			\begin{enumerate}
				\item Compute $T_{2_i}^{(m)}$, using updated values $U_j^{(m)}$ for $j<i$ and $U_j^{(m-1)}$ for $j>i$ \hfill (Eq.~\eqref{eq:T2})
				\item Compute $U_{1_i}$ and $U_{2_i}$ \hfill (Eqs.~\eqref{eq:prop11}--\eqref{eq:prop22})
				\item Compute $\alpha_i, \kappa_i, \lambda_i$ \hfill (coefficients of the quadratic equation)
				\item Compute $U_{3_i}$ \hfill (Eq.~\eqref{eq:prop33})
				\item \textbf{if} $U_{1_i}+R(\bar p_{\mathrm{ext}}) \le \bar u_v$ \textbf{then} $U_i^{(m)} \gets U_{1_i}$
				\item \textbf{else if} $U_{2_i}+R(\bar p_{\mathrm{ext}}) \ge \bar u_c$ \textbf{then} $U_i^{(m)} \gets \omega\,U_{2_i} + (1-\omega)\,U_i^{(m-1)}$
				\item \textbf{else} $U_i^{(m)} \gets U_{3_i}$
			\end{enumerate}
			\item \textbf{end for}
			\item Compute the relative error $e^{(m)}=\|U^{(m)}-U^{(m-1)}\|/\|U^{(m)}\|$
			\item $m \gets m+1$
		\end{enumerate}
		\item \textbf{until} $e^{(m)}<\varepsilon$
		\item \textbf{return} $U^{(m)}$
	\end{enumerate}
	\caption{Gauss--Seidel algorithm for solving the nonlinear system (NGS)}
	\label{alg:gs_u}
	\addtocounter{figure}{-1}
\end{figure}

\begin{figure}[H]
	\renewcommand{\figurename}{Algorithm}
	\stepcounter{algo}
	\renewcommand{\thefigure}{\arabic{algo}}
	\renewcommand{\theHfigure}{alg.\arabic{algo}}
	\small
	\begin{enumerate}
		\item[] \textbf{Require:} initial guess $U^{(0)}$, tolerance $\varepsilon_{\mathrm{tol}}$, parameters $\omega_{\max}$, $C_\omega$
		\item[] \textbf{Ensure:} approximate solution $U^{(m)}$ and pressure $\bar p=R^{-1}\bigl(U^{(m)}\bigr)$
		\item Assemble $M$ and $E$
		\item $m\gets 0$
		\item \textbf{while} $\|\mathcal R(U^{(m)})\|_2>\varepsilon_{\mathrm{tol}}$ \textbf{do}
		\begin{enumerate}
			\item $\mathcal R(U^{(m)})\gets MU^{(m)}-E\rho(U^{(m)})$
			\item $J(U^{(m)})\gets M-E\,D_{\rho'}(U^{(m)})$ \hfill (Eq.~\eqref{eq:newton_jacobian})
			\item $\omega_m\gets\min\bigl(\omega_{\max},\,C_\omega/\|\mathcal R(U^{(m)})\|_2\bigr)$
			\item Solve the sparse linear system $J(U^{(m)})\,\delta U^{(m)}=-\mathcal R(U^{(m)})$
			\item $U^{(m+1)}\gets U^{(m)}+\omega_m\,\delta U^{(m)}$
			\item $m\gets m+1$
		\end{enumerate}
		\item \textbf{end while}
	\end{enumerate}
	\caption{Newton--Raphson method with adaptive relaxation}
	\label{alg:newton_adaptive}
	\addtocounter{figure}{-1}
\end{figure}


\end{document}